%
\magnification=\magstep1   
\input amstex
\UseAMSsymbols
\input pictex
\vsize=23truecm
\NoBlackBoxes
\parindent=18pt
  
   \font\rmk=cmr8    \font\itk=cmti8  \font\ttk=cmtt8

\def\op{\text{\rm op}}
\def\mod{\operatorname{mod}}

\def\Hom{\operatorname{Hom}}

\def\Ext{\operatorname{Ext}}

\def\rad{\operatorname{rad}}

\def\Ker{\operatorname{Ker}}
\def\Cok{\operatorname{Cok}}

\def\Im{\operatorname{Im}}
\def\diff{\operatorname{diff}}
\def\perf{\operatorname{perf}}
\def\bdim{\operatorname{\bold{dim}}}
  
\def\arr#1#2{\arrow <1.5mm> [0.25,0.75] from #1 to #2}
	\vglue1truecm

\plainfootnote{}
{\rmk 2010 \itk Mathematics Subject Classification. \rmk 
Primary 16G10, 16G50, 16E65.  Secondary 16G70, 18G25.}

\centerline{\bf Representations of quivers over the algebra of dual numbers.}
                     \bigskip
\centerline{Claus Michael Ringel and Pu Zhang}     
		  \bigskip\medskip

{\narrower\narrower
Abstract: 
The representations of a quiver 
$Q$ over a field 
$k$ (the $kQ$-modules, where $kQ$ is the path algebra of $Q$
over $k$)  have been studied for a long time, 
and one knows quite well the structure of the module category $\mod kQ$. It seems to be worthwhile to consider also representations of $Q$ over arbitrary finite-dimensional $k$-algebras $A$.
Here we draw the attention to the case when $A = k[\epsilon]$ is the algebra
of dual numbers (the factor algebra of the polynomial ring $k[T]$ in one variable $T$ modulo
the ideal generated by $T^2$), thus to the $\Lambda$-modules, where $\Lambda = kQ[\epsilon] =
kQ[T]/\langle T^2\rangle.$ The algebra $\Lambda$ is a $1$-Gorenstein algebra, thus the torsionless 
$\Lambda$-modules are known to be of special interest (as the Gorenstein-projective or
maximal Cohen-Macaulay modules). They form a Frobenius category
$\Cal L$, thus the corresponding stable category $\underline{\Cal L}$
is a triangulated category. As we will see, the category $\Cal L$ is the category of
perfect differential $kQ$-modules and $\underline{\Cal L}$ is the corresponding homotopy
category. The category $\underline{\Cal L}$
is triangle equivalent to
the orbit category of the derived category $D^b(\mod kQ)$ modulo the shift and
the homology functor 
$H\: \mod \Lambda \to \mod kQ$ yields a bijection between the
indecomposables in $\underline{\Cal L}$ and those in $\mod kQ$. Our main interest lies
in the inverse, it is given by the minimal $\Cal L$-approximation. Also, we 
will determine the kernel of the restriction of the functor $H$ to $\Cal L$ and 
describe 
the Auslander-Reiten quivers  of $\Cal L$ and $\underline{\Cal L}.$
\par} 
	\bigskip
Throughout the paper, $k$ will be a field and $Q$ will be a finite connected directed quiver. 
The starting point for the considerations of this paper is the following result which concerns 
the structure of the homotopy category of perfect differential $kQ$-modules.  
This assertion should be well-known, but we could not find a reference.

Let us recall that given a ring $R$, a differential $R$-module is by definition a pair
$(N,\epsilon)$ where $N$ is an $R$-module and $\epsilon$ an endomorphism of $N$ such that 
$\epsilon^2 = 0.$
If $(N,\epsilon)$ and $(N',\epsilon')$ are differential $R$-modules, 
a morphism $f\:(N,\epsilon) \to (N',\epsilon')$ is given 
by an $R$-linear map $f\:N \to N'$ such that $\epsilon' f = f\epsilon$. 
The morphism $f\:(N,\epsilon) \to (N',\epsilon')$ is said to be {\it homotopic to zero}
provided there exists an $R$-linear map $h\:N \to N'$ such that $f = h\epsilon + \epsilon'h.$
A differential
$R$-module $(N,\epsilon)$ is said to be {\it perfect} provided $N$ is a finitely generated
projective $R$-module.
We denote by ${\diff}_{\perf}(R)$ the category of perfect differential $R$-modules, and by
$\underline{\diff}_{\perf}(R)$ the corresponding homotopy category.
Let us denote by $H$ the homology functor: it attaches to a differential $R$-module $(N,\epsilon)$
the $R$-module $H(N,\epsilon) = \Ker \epsilon/\Im \epsilon.$ It is well-known that
$H$ vanishes on the maps which are homotopic to zero. 

If $R$ is noetherian, let us denote by $D^b(\mod R)$ the bounded derived category 
of finitely generated $R$-modules. This is a triangulated category and its shift functor
will be denoted by $[1]$. 
	\medskip
{\bf Theorem 1.} (a) {\it The category ${\diff}_{\perf}(kQ)$ 
of perfect differential $kQ$-module is a Frobenius category
whose stable category $\underline{\diff}_{\perf}(R)$ 
is the orbit category $D^b(\mod kQ)/[1]$.}

(b) {\it The homology functor $H\: {\diff}_{\perf}(kQ) \to \mod kQ$ 
is a full and dense functor which furnishes a bijection between the 
indecomposables in the homotopy category $\underline{\diff}_{\perf}(kQ)$ and those in $\mod kQ$. 
It yields a quiver embedding $\iota$ of the
Auslander-Reiten quiver of $\mod kQ$ into the Auslander-Reiten quiver of 
the homotopy category $\underline{\diff}_{\perf}(kQ)$.}
		\medskip
We should remark that the study of differential modules themselves may have been neglected
by the algebraists, however it is clear that the graded version, namely complexes, play 
an important role in many 
parts of mathematics. Theorem 1 is an immediate consequence of 
well-known results concerning perfect complexes over $kQ$: the category of
perfect complexes is a Frobenius category, thus the corresponding stable category
is $D^b(\mod kQ)$; the homology functor $H_0$ from the
category of perfect complexes to $\mod kQ$ is full and dense and it furnishes a bijection between the 
shift orbits of the indecomposables in the homotopy category of perfect complexes  
and the indecomposables in $\mod kQ$; also, it yields an embedding $\iota$ of the
Auslander-Reiten quiver of $\mod kQ$ into the Auslander-Reiten quiver of 
the homotopy category of the perfect complexes.  

Theorem 1 follows from these assertions, 
using the covering theory as
developed by Gabriel and his school (or the equivalent theory of
group graded algebras by Gordon and Green); namely one just looks at the 
forgetful functor
from the category of perfect complexes of $kQ$-modules to the
category of perfect differential $kQ$-modules.  
The main property to be used
is the local boundedness of the category of perfect complexes. For the fact that
the orbit category $D^b(\mod kQ)/[1]$ is triangulated, we may refer to the general criterion
given recently by Keller [Ke], however one also may use directly the Frobenius category structure
of the category of perfect differential $kQ$-modules.  
Some further comments on this proof will be given in section 3. 
	\medskip
The proper framework for  Theorem 1 seems to be the category
of {\bf all} differential $kQ$-modules, this category may be interpreted in several ways.
It is the category of $\Lambda$-modules,
where $\Lambda = kQ[\epsilon] = kQ[T]/\langle T^2\rangle$ (with $T$ a central variable). Note
that  $\Lambda = AQ$
may be considered as the path algebra of the quiver $Q$ 
over the 2-dimensional local algebra 
$A = k[\epsilon] = k[T]/\langle T^2\rangle$ of dual numbers, thus the $\Lambda$-modules are 
just the representations of
$Q$ over the ring $A$. Also, we may write $\Lambda$ as the tensor product
of $kQ$ with $A$ over $k$. 

The aim of this paper to analyze Theorem 1 as dealing with two 
subcategories of the module category $\mod \Lambda$, the subcategories
of interest are on the one hand the category $\mod kQ$ (these are the $\Lambda$-modules 
annihilated by $\epsilon$), and the category of perfect differential $kQ$-modules on the 
other hand. 

The decisive property which we will use is the fact that $\Lambda$ is a $1$-Gorenstein algebra. For any $1$-Gorenstein algebra $\Lambda$, 
the category $\Cal L = \Cal L_\Lambda$ of the torsionless $\Lambda$-modules
is of interest,
these are the $\Lambda$-modules which are submodules of projective modules,
but they also can be characterized differently: 
the torsionless $\Lambda$-modules are just the 
Gorenstein-projective or maximal Cohen-Macaulay modules as considered in [EJ,B]),
and the modules of $G$-dimension 0 in the sense [AB].
Note that $\Cal L$ is a Frobenius category, thus the corresponding stable category 
$\underline{\Cal L}$ (obtained from $\Cal L$ by factoring out all the maps which factor
through the subcategory $\Cal P$ of the projective $\Lambda$-modules) is a triangulated
category. In our case $\Lambda = kQ[\epsilon]$, the
category $\Cal L$ is precisely the category of perfect differential $kQ$-modules,
$$
 \Cal L = {\diff}_{\perf}(kQ)\quad \text{and}\quad 
 \underline{\Cal L} = 
\underline{\diff}_{\perf}(R),
$$
and every module in $\Cal L$ is even strongly Gorenstein-projective, see 4.10.

The basic functor to be considered is the homology functor $H\: \mod \Lambda \to \mod kQ,$ it
sends a representation $M = (M_i,M_\alpha)_{i\in Q_0,\alpha\in Q_1}$ 
of the quiver $Q$ over $k[\epsilon]$ to the homology with respect
to the action of $\epsilon$, thus, $H(M) = (H(M_i),H(M_\alpha))_{i,\alpha}$.
Besides the functor $H\:\mod \Lambda \to \mod kQ$ we will consider a reverse 
construction $\eta$ which is not functorial (but of
course stably functorial), the minimal right $\Cal L$-approximation.
	\medskip
{\bf Theorem 2.} {\it The algebra $\Lambda = kQ[\epsilon]$ is $1$-Gorenstein. The functor 
$H:\Cal L \to \mod kQ$ is full and 
induces a bijection between the indecomposable $\Lambda$-modules
in $\Cal L \setminus \Cal P$ and the indecomposable $kQ$-modules.
The inverse bijection is given by taking the minimal right $\Cal L$--approximation
of an indecomposable $kQ$-module.}
	\bigskip
We obtain an embedding $\iota$ of the
Auslander-Reiten quiver $\Gamma(\mod kQ)$ 
of $\mod kQ$ into the Auslander-Reiten quiver $\Gamma(\underline{\Cal L})$ 
of $\underline{\Cal L}$
by sending an indecomposable $kQ$-module $N$ to $\eta N$. The only arrows which are not
obtained in this may are the following ``ghost maps'':
For any vertex $y$ of $Q$, let $P_0(y)$ and $I_0(y)$ be the
corresponding indecomposable projective or injective $kQ$-module, respectively. 
For any  $y$, we construct a homomorphism
$c(y)\:\eta(I_0(y)/S(y)) \to P_0(y)$. These homomorphisms 
yield the arrows in $\Gamma(\underline{\Cal L})$ 
which are not in the image of $\iota$. 
In addition, we show that $\tau_{\Cal L} P_0(y) = \eta I_0(y)$, where $\tau_{\Cal L}$ is 
the Auslander-Reiten translation $\tau_{\Cal L}$ of $\Cal L$, 
In this way we get the required extension of the translation map of $\Gamma(\mod kQ)$. 

Of course, in order to obtain
the Auslander-Reiten quiver of $\Cal L$ itself, we have to add the indecomposable
projective $\Lambda$-modules $P(y) = P_0(y)\otimes_k A$. 
But this is easy: $\rad P(y)$ belongs to $\Cal L$ and is
indecomposable, actually $H(\rad P(y))$ is  
just the simple module $S(y)$ corresponding to the vertex $y$. 
	\medskip
{\bf Theorem 3.} {\it The kernel of the functor $H\:\Cal L \to \mod kQ$ is a
finitely generated ideal of $\Cal L$, it is generated by the identity morphisms of the indecomposable projective $\Lambda$-modules and the homomorphisms 
$\eta I(x) \to P(y)$, where $x,y$ are vertices of $Q$.}
	\medskip
In fact, instead of using all the maps $\eta I(x) \to P(y)$, it is sufficient to take the
maps $c(y)$ mentioned above.  
	\bigskip
Here is an outline of the paper. The first section describes the context
of this investigation. In section 2, we show that the perfect differential $kQ$-modules
are precisely the torsionless $kQ[\epsilon]$-modules, thus the Gorenstein-projective
$kQ[\epsilon]$-modules. Section 3 provides some details for the covering approach. Section 4
is the central part, here we discuss in which way the homology functor $H$ and the 
$\Cal L$-approximation $\eta$ 
are inverse to each other. Sections 5 and 6 deal with the ghost maps, section 7
with the position of the indecomposable projective $\Lambda$-modules in the 
Auslander-Reiten sequences of $\Cal L$.
Sections 8 to 10 are devoted to examples and further remarks.
	\bigskip\bigskip
\vfill\eject
{\bf 1. The context.}
	\medskip
{\bf 1.1.} An explicit description of the category of Gorenstein-projective modules
is known only for very few algebras. In a recent paper Luo and Zhang [LZ] gave 
a characterization of the 
Gorenstein-projective $AQ$-modules, where $Q$ is a finite directed quiver and $A$
a any $k$-algebra, thus one may try to use this result in order to construct these
modules explicitly. The present paper deals with the very special case of the algebra
$A = k[\epsilon]$ of dual numbers, this may be considered as an interesting test case.

Let us stress that the class of $1$-Gorenstein algebras is a class of
algebras which includes both the hereditary and the self-injective algebras --- two classes
of algebras whose representations have been investigated 
very thoroughly and have been shown to be strongly related to Lie theory. 
Thus one might hope that all the $1$-Gorenstein algebras have such a property.
Now recently,  Keller and Reiten [KR] identified another class of $1$-Gorenstein
algebras, namely the cluster tilted algebras, and this again is a class 
of algebras related to Lie theory. Of course, the result of the present paper
also supports the hope.
	
As we have mentioned, an explicit description of the category of Gorenstein-projective modules
is known only in few cases. Chen ([C], see also [RX]) 
recently has shown that for $\Lambda$ of Loewy length
at most 2, the stable category of Gorenstein-projective $\Lambda$-modules is a union of
categories with Auslander-Reiten quiver of tree type $\Bbb A_1$, 
thus not very exciting. Now, the next case of interest are the
artin algebras of Loewy length 3, and for this case we provide a wealth of examples. 
Namely, if $Q$
is a finite bipartite quiver (i.e., all vertices are sinks or sources), then the
algebra $kQ[\epsilon]$ is of Loewy length at most $3$.  
	\medskip
{\bf 1.2.} 
The theorems allow to transfer a lot of results known for the module category $\mod kQ$ 
to the category $\Cal L$ of Gorenstein-projective $\Lambda$-modules, where 
$\Lambda = kQ[\epsilon].$ For example, the Kac Theorem [K] yields:
	\medskip
{\it The homology dimension vector $\bdim H(-)$ maps any indecomposable object
in $\Cal L \setminus \Cal P$ to a positive root of the corresponding Kac-Moody algebra $\bold g$.
For any positive real root $\bold r$, there is a unique isomorphism class of indecomposable
modules $M$ in $\Cal L$ with
$\bdim H(M) = \bold r$; if $k$ is an infinite field, then for every positive imaginary root 
$\bold r$ of 
$\bold g$, there are infinitely many isomorphism classes of indecomposable
modules $M$ in $\Cal L$ with
$\bdim H(M) = \bold r$.}
	\medskip
{\bf 1.3.} 
The relationship between abelian categories and triangulated categories
has always been considered as fascinating, but also mysterious. It was clear from
the beginning that starting with a suitable abelian category $\Cal A$ (namely
the module category $\Cal A$ of a self-injective algebra)  one may obtain
a triangulated categories by factoring out some finitely generated ideal
(namely the ideal
of all maps which factor through a projective module). 
Only quite recently, it was observed
that there are also examples in the reverse direction: if one starts with a cluster
category (a triangulated category, according to [Ke]) and factors out the ideal of all maps which
factor through the additive category generated by a fixed cluster-tilting object, 
then one obtains an abelian category, namely
the module category of the corresponding cluster-tilted algebra (thus
an abelian category), see Buan-Marsh-Reineke-Reiten-Todorov [BMRRT].

The results of the present paper should be seen in this context. As in the case of the cluster categories we start with a triangulated category $\Cal T$ and factor out a finitely
generated ideal of $\Cal T$ in order to obtain an abelian category. But whereas in the
cluster category case the ideal is generated by some identity maps, here we deal with an
ideal $I$ which lies inside the radical of $\Cal T.$ According to 6.5, we even
know that $I^2 = 0.$
	
It seem to be of interest that actually we deal with two related
subcategories of a module category, one is a Frobenius category $\Cal F$, the other an abelian
category $\Cal A$, such that there is a canonical bijection between the indecomposable objects
in the stable category $\underline{\Cal F}$ and the indecomposable objects in $\Cal A.$
Of course, if we consider for a self-injective algebra $R$ 
the subcategories $\Cal F = \mod R$ 
and $\Cal A = \mod R/S(y) R$, then there is such a canonical bijection, namely the identity.
The examples which we consider are more intricate: here, 
$\Cal F = \Cal L$ is the category of Gorenstein-projective $kQ[\epsilon]$-modules, 
$\Cal A$ is the category of $kQ$-modules,
and the bijection is given by the functor $H$ and the $\Cal L$-approximation $\eta$.
	\medskip
{\bf 1.4.} A long time ago,
it has been shown by Buchweitz [B] that given a Gorenstein algebra $\Lambda$, 
the Verdier quotient of the bounded derived category $D^b(\mod\Lambda)$ modulo
the subcategory of perfect complexes can be identified with the 
stable category of Gorenstein-projective $\Lambda$-modules, and Orlov [O] proposed the name
{\it triangulated category of singularities} for this Verdier quotient.
In our case $\Lambda = kQ[\epsilon]$, we show that the triangulated category of singularities
is just $D^b(\mod kQ)/[1]$.
		\bigskip\bigskip
{\bf 2. The basic observation.} 
	\medskip	
Let $R$ be an arbitrary ring. As above, we define $R[\epsilon] = R[T]/\langle T^2\rangle,$ where $T$
is a variable which is supposed to commutate with all the elements of $R$. The
$R[\epsilon]$-modules are just the {\it differential} $R$-modules, 
they may be written as $(N,f),$
where $N$ is an $R$-module and $f\:N \to N$ is an $R$-endomorphism with $f^2 = 0$ (namely,
if such a pair $(N,f)$ is given, then $N$ can be considered as an $R[\epsilon]$-module by
defining the action of $\epsilon$ on $N$ as being given by $f$); by abuse of notation, we sometimes
will write $\epsilon$ instead of $f$. A differential $R$-module $(N,f)$
will be said to be {\it perfect} provided $N$ is a finitely generated projective $R$-module.

For  an $R$-module $N$, let $N[\epsilon] =
(N\oplus N; \left[\smallmatrix 0 & 1 \cr 0 & 0 \endsmallmatrix\right])$, this is an $R[\epsilon]$-module
(note that if we take $N = R$, then the module $R[\epsilon]$ is just the regular
representation of the ring $R[\epsilon]$, thus there is no conflict of notation).
If $N$ is a finitely generated
$R$-module, then $N[\epsilon]$ is a finitely generated $R[\epsilon]$-module; 
if $N$ is a projective
$R$-module, then $N[\epsilon]$ is a projective $R[\epsilon]$-module. 
In particular, if $N$ is finitely generated and projective, then $N[\epsilon]$ is a perfect 
$R[\epsilon]$-module. But a perfect $R[\epsilon]$-module may not be of the form $N[\epsilon]$.
Also note that a finitely generated projective $R[\epsilon]$ module is perfect, but the converse
is not true. 

Let us recall that an $R$-module is said to be {\it torsionless} provided it is a submodule
of a projective $R$-module. Thus, a differential $R$-module is torsionless
if it is a submodule of a projective $R[\epsilon]$-module.  
Also recall that a ring $R$ is {\it left hereditary} provided any
torsionless left module is projective.
	\medskip
{\bf 2.1. Lemma.} {\it Let $R$ be left noetherian and left hereditary.
A differential $R$-module is finitely generated and torsionless
if and only if it is perfect.}
	\medskip
Proof. Let $(N,f)$ be a differential $R$-module. 

First, assume that $(N,f)$ is finitely generated and torsionless. Since $(N,f)$ is torsionless,
we know that $(N,f)$ is a submodule of a free $R[\epsilon]$-module, and since $(N,f)$
is finitely generated, it is even a submodule of a free $R[\epsilon]$-module of finite rank.
Thus $(N,f)$ can be embedded into a module of the form $(R^t\oplus R^t,
\left[\smallmatrix 0 & 1 \cr
                     0 & 0 \endsmallmatrix\right])$ for some natural number $t$.
In particular, we see that $N$ is a submodule of $R^{2t}$. Since $R$ is left noetherian
and hereditary, we conclude that $N$ is a finitely generated projective $R$-module.

Conversely, assume that $N$ is a finitely generated projective $R$-module.
We denote by $N'$ the kernel, by $N''$ the image of $f$; let $u'\:N'' \to N'$ and 
$u\:N' \to N$ be the inclusion maps. Note that $f$ induces an epimorphism $f'\:N \to N''$
with $f = uu'f'.$ Since $R$ is left hereditary and $N$ is a projective $R$-module, we see
that $N''$ is also projective. It follows that there is a homomorphism $s\:N'' \to N$ such that
$f's = 1_{N''}$, and, as a consequence the map 
$$
 \bmatrix u & s\endbmatrix\:N'\oplus N'' \to N
$$
is an isomorphism. But $uu' = uu'f's = fs$ shows that
$$
 \bmatrix u & s\endbmatrix\bmatrix 0 & u'\cr 0 & 0 \endbmatrix 
  = f \bmatrix u & s\endbmatrix,
$$
therefore $\bmatrix u & s\endbmatrix$ is an
isomorphism between $(N'\oplus N'',\left[
\smallmatrix 0 & u'\cr 0 & 0\endsmallmatrix\right])$ and $(N,f)$. It 
remains to observe that there is the following embedding of differential $R$-modules:
$$
 \bmatrix 1 & 0 \cr 0 & u'\endbmatrix\:(N'\oplus N'',\bmatrix 0 & u'\cr 0 & 0 \endbmatrix)
 \to (N' \oplus N',\bmatrix 0 & 1\cr 0 & 0 \endbmatrix),
$$
and that 
$(N' \oplus N',\left[\smallmatrix 0 & 1\cr 0 & 0 \endsmallmatrix\right])$ is a projective $\Lambda$-module.
This shows that 
$(N'\oplus N'',\left[\smallmatrix 0 & u'\cr 0 & 0 \endsmallmatrix\right])$ and therefore
$(N,f)$ is torsionless. Since $N$ is a finitely generated $R$-module, we see that $(N,f)$
is a finitely generated $\Lambda$-module. \hfill $\square$
	\medskip
Let us consider now finite-dimensional $k$-algebras, where $k$ is a field.
Recall [AR1] 
that such an algebra $R$ is called a {\it Gorenstein algebra of Gorenstein dimension $1$}
or just a {\it $1$-Gorenstein algebra}
provided the injective dimension of ${}_RR$ as well as $R_R$ is equal to $1$.
Since $\Lambda$ is the tensor product of
a hereditary and a self-injective algebra, one sees immediately that $\Lambda$
is $1$-Gorenstein.

Given any finite-dimensional $k$-algebra $\Lambda$, a $\Lambda$-module $M$ is said to be
{\it Gorenstein-projective} [EJ] provided there exists an exact (not necessarily bounded) complex 
$P_\bullet = (P_i,\delta_i)_i$ of finitely generated projective $\Lambda$-modules
such that also $\Hom_\Lambda(P_\bullet,{}_\Lambda\Lambda)$ is exact and such that 
$M$ is the image of $\delta_0.$ If $\Lambda$ is Gorenstein, then a $\Lambda$-module $M$
is Gorenstein-projective if and only if $\Ext^i(M,\Lambda) = 0$ for all $i\ge 1.$

The following proposition is well-known.
	\medskip
{\bf 2.2. Proposition.} {\it Let $R$ be a finite-dimensional $k$-algebra which is hereditary.
Then $\Lambda = R[\epsilon]$  is a Gorenstein algebra of Gorenstein dimension $1$. 
If $M$ is any $\Lambda$-module. then the following conditions are equivalent:
\item{\rm (i)} $M$ is Gorenstein-projective,
\item{\rm (ii)} $M$ is torsionless,
\item{\rm (iii)} $\Ext^1_\Lambda(M,\Lambda) = 0.$}
	\medskip
{\bf 2.3.} 
{\it Given any $\Lambda$-module $M$, there is an exact sequence 
$$
 0 @>>> P @>>> \eta M @>g>> M @>>> 0,
$$
such that $P$ is projective, $\eta M$ belongs to $\Cal L$ and $g$ is a right minimal map}
(see for example [EJ], Theorem 11.5.1). 

	\medskip
The map $g$ (or also the module $\eta M$) is called the minimal right $\Cal L$-approximation,
since any map $h\: L \to M$ with $L\in \Cal L$, can be factorized as $h = h'g$ with
$h'\:L \to \eta M$. Of course, this factorization property and the minimality of $g$
implies that $\eta M$ is uniquely determined by $M$, up to an isomorphism, but there
there is not necessarily a canonical isomorphism. Also, one should note that $\eta M$
is the universal extension of $M$ from below, using projective modules.
	\bigskip
Section 4 below will be devoted to a detailed study of the minimal
right $\Cal L$-approximations $\eta N \to N$, where $\Cal L$ is the category
of torsionless $\Lambda$-modules with $\Lambda = kQ[\epsilon]$, and $N$
is a $\Lambda$-module which is annihilated by $\epsilon$ (thus a $kQ$-module).
	\bigskip

{\bf 2.4.} Again, let $\Lambda = kQ[\epsilon]$ and $\Cal L$ 
the category of torsionless $\Lambda$-modules.
The Frobenius structure on $\Cal L$ is given by those sequences
$\eta = (0 @>>> X @>f>> Y @>g>> Z @>>> 0)$ in $\Cal L$ 
which yield exact sequences when we apply
$\Hom(\Lambda,-)$ and $\Hom(-,\Lambda)$. To say that we obtain an exact sequence
when we apply $\Hom(\Lambda,-)$ means that the given sequence is exact. 
And to say that we obtain an exact sequence when we apply $\Hom(-,\Lambda)$
means that $\Hom(f,\Lambda)$ is surjective, thus that $f$ is a left
$\Lambda$-approximation. But if $\eta$ is any exact sequence in $\mod \Lambda$
with $Z$ in $\Cal L$, then $f$ is a left $\Lambda$-approximation, since
$\Ext^1(L,\Lambda) = 0.$ This shows that the Frobenius structure on $\Cal L$
is given just by the exact sequences of $\mod \Lambda$ with $X,Y,Z \in \Cal L.$.
	
Given an exact sequence $\eta = (0 @>>> X @>f>> Y @>g>> Z @>>> 0)$ in $\mod \Lambda$
with $Z \in \Cal L,$ then $\eta$ considered as an exact sequence of $kQ$-modules
splits (since $Z\in \Cal L$ means that $Z$ is a perfect differential $kQ$-module,
thus that $Z$ is $kQ$-projective).
	\bigskip
{\it The functor $H\:\underline{\Cal L} \to \mod kQ$ is a cohomological functor.}
	\medskip
Proof: We have to show that any triangle $X @>>> Y @>>> Z @>>>$ in
$\underline{\Cal L}$ yields under $H$ an exact sequence $H(X)  \to H(Y) \to H(Z).$
Recall that such a triangle in $\underline{\Cal L}$ starting with 
a homomorphism $f\:X \to Y$ in $\Cal L$ can be  constructed as follows: 
$$
\CD
 0 @>>> X @>u>> P @>>> \Sigma X @>>> 0 \cr
 @.    @VfVV    @Vf'VV    @| \cr
 0 @>>> Y @>g>> Z @>>> \Sigma X @>>> 0,
\endCD
$$
this is a commutative diagram in $\mod\Lambda$ with exact rows such that 
$u$ is a left $\Lambda$-approximation (therefore, $P$ is a projective $\Lambda$-module and $\Sigma X \in \Cal L$).
Thus $Z = (Y\oplus P)/U$, where $U = \{(f(x),-u(x))\mid x\in X\}$ and
$g(y) = \overline{(y,0)}),$
$f'(p) = \overline{(0,p)}),$ for $y\in Y$ and $p\in P.$ 
The exact sequence 
$$
 0 @>>> X @>{\left[\smallmatrix f \cr u\endsmallmatrix\right]}>> 
       Y\oplus P @>{\left[\smallmatrix g & f'\endsmallmatrix\right]}>> Z @>>> 0
$$
yields under $H$ a sequence
$$
 H(X) @>{\left[\smallmatrix H(f) \cr H(u)\endsmallmatrix\right]}>> 
       H(Y)\oplus H(P) @>{\left[\smallmatrix H(g) & H(f')\endsmallmatrix\right]}>> H(Z)
$$
with zero composition. But $H(P) = 0$, thus $H(g)H(f) = 0.$ It remains to
be seen that any element of  the kernel of $H(g)$ is in the image of $H(f)$.
An element of $H(Y)$ is the residue class $\overline y$ modulo $\epsilon Y$
of an element $y\in Y$ with $\epsilon y = 0$, and $\overline y$ belongs to the
kernel of $H(g)$ iff $g(y)$ belongs to $\epsilon Z.$ Thus assume that
$g(y) = \epsilon z$ with $z \in Z$, say $z = (y',p)+U$ for some $y'\in Y$
and $p\in P.$ Thus $(y,0)-(\epsilon y',\epsilon p)$ belongs to $U$, this means that
there is $x\in X$ with $y-\epsilon y' = f(x)$ and $0-\epsilon p = -u(x)$.
Note that $u(\epsilon x) = \epsilon u(x) = \epsilon^2(p) = 0.$ Since $u$
is injective, $\epsilon x = 0$, thus $x$ belongs to the kernel of $\epsilon.$
The equality $y = f(x)+\epsilon y'$ shows that $\overline y =
\overline{f(x)} = H(f)(\overline x)$ in
$H(Y)$, thus $\overline y$ is in the image
of $H(f)$. \hfill $\square$

	\bigskip\bigskip

{\bf 3. Proof of theorem 1.}
	\medskip
Given any ring $R$, we denote by $\Cal P = \Cal P_R$ the category of finitely generated projective
$R$-modules. 

We are interested in complexes of $R$-modules, such complexes may be considered as
differential graded $R$-modules. In particular, we will consider perfect complexes
(or perfect differential graded $R$-modules), these are the bounded complexes which use
only finitely generated projective $R$-modules.
We denote by $C^b(\Cal P_R)$ the category of perfect complexes, and by $K^b(\Cal P_R)$
the corresponding homotopy category. Let us stress that $C^b(\Cal P_R)$ is a Frobenius
category and $K^b(\Cal P_R)$ is just the corresponding stable category,
say with stabilization functor $\pi\:C^b(\Cal P_R) \to K^b(\Cal P_R)$, it
sends a map to its homotopy class.

Note that in case $R = kQ,$ where $Q$ is a finite directed
quiver, the ring $R$ has finite global dimension, and
$K^b( \Cal P_R)$ can be identified with the bounded derived category
$D^b(\mod R).$ 
In dealing with categories of complexes, the shift functor will be denoted by $[1]$.

There is the following commutative diagram
$$
\CD
 C^b(\Cal P_{R}) @>\pi>>  K^b(\Cal P_{R})\cr
 @V\gamma VV                      @VV\gamma V  \cr
{\diff}_{\perf}(R)
  @>\pi>> \underline{\diff}_{\perf}(R)
\endCD
$$ 
here, the horizontal functors $\pi$ are just the stabilization functors, whereas the
vertical functors $\gamma$ are obtained by forgetting the grading (such a forgetful
functor is sometimes called a compression functor or,
in the covering theory of the Gabriel school, the corresponding pushdown functor).

What is important here, is the fact that for $R = kQ$
the category $C^b(\Cal P_{R})$ is locally bounded,
thus the functors $\gamma$ are dense and provide a bijection between the
shift-orbits of the indecomposable objects in $C^b(\Cal P_{R})$ and the
indecomposable objects in ${\diff}_{\perf}(R)$
and similarly, between the 
shift-orbits of the indecomposable objects in $K^b(\Cal P_{R})$ and the
indecomposable objects in $\underline{\diff}_{\perf}(R)$.

Let us denote by $H_0 \: C^b(\Cal P_{R}) \to \mod R$ the functor which attaches
to a complex $P_\bullet = (P_i,\delta_i)$ the homology at the position $0$.
Then it is well-known [H1] that $H_0$ provides a bijection between the shift orbits
of indecomposable objects in $K^b(\Cal P_{kQ})$ and the indecomposable $kQ$-modules.
This completes the proof.  \hfill $\square$
	\bigskip
For the case $R = kQ$, let us exhibit the diagram above again, but now using the notation 
$\Cal L = {\diff}_{\perf}(kQ)$ and $\underline{\Cal L} = \underline{\diff}_{\perf}(kQ)$
which is used in the further parts of the paper:
$$
\CD
 C^b(\Cal P_{kQ}) @>\pi>>  K^b(\Cal P_{kQ})\cr
 @V\gamma VV                      @VV\gamma V  \cr
 \Cal L
  @>\pi>> \underline{\Cal L}
\endCD
$$ 

 	\bigskip\bigskip
{\bf 4. Explicit bijections.}
	\medskip
Let $Q$ be a finite directed quiver, $kQ$ its path algebra and $\Lambda = kQ[\epsilon].$
Recall that we denote by $\Cal L$ the category of  (finitely generated)
torsionless modules. 
Given any $\Lambda$-module $M$, let $\eta M \to M$ be a minimal right $\Cal L$-approximation.
	\bigskip
We consider $\mod kQ$ as the full subcategory of $\mod \Lambda$ given by all the
$\Lambda$-modules $N$ which are annihilated by $\epsilon.$
	\bigskip
{\bf 4.1.} {\it The $\Lambda$-modules in $\Cal L \cap \mod kQ$ are just the projective $kQ$-modules.}
	\medskip
Proof.  Let $M$ be a module in $\Cal L \cap \mod kQ$. Since $M$ is in $\Cal L$,
we may consider $M$ as a submodule of $({}_\Lambda\Lambda)^t$ for some $t$. Since $M$ is
annihilated by $\epsilon$, it must lie inside the submodule 
$(\epsilon\Lambda)^t$ (since $\{x\in \Lambda\mid
\epsilon x = 0\} = \epsilon \Lambda$). But $\epsilon \Lambda \simeq kQ$ as $\Lambda$-modules
and thus also as $kQ$-modules.
Now 
$M$ is a submodule of $kQ$, thus it is a projective $kQ$-module. Of course, conversely, a projective
$kQ$-module belongs to $\Cal L$.  \hfill $\square$
	\medskip
{\bf 4.2.} {\it Let $M$ be in $\Cal L$, let
$M' = \{x\in M\mid \epsilon x = 0\}$ and $M'' = \epsilon M.$ Then $M'' \subseteq M'$,
$H(M) = M'/M''$ and the exact sequence
$$
 0 \to M'' \to M' \to H(M) \to 0
$$
is a projective $kQ$-resolution of $H(M)$.}  
	\medskip
Proof.  Since $M'$ is a submodule of $M$, it belongs to $\Cal L.$ Since $M'$ is annihilated
by $\epsilon,$ is belongs to $\mod kQ$, thus $M'$ is a projective $kQ$-module by 4.1.
 \hfill $\square$
	\medskip

{\bf 4.3.} {\it Let $M$ be a $\Lambda$-module. If $M$ is projective, then
$H(M) = 0$ and $\epsilon M$ is a projective $kQ$-module. Conversely, if 
$H(M) = 0$ and $N = \epsilon M$ is a projective
$kQ$-module, then $M$ is isomorphic to $N[\epsilon]$, thus $M$ is a projective $\Lambda$-module.}
	\medskip
Proof. If $M = {}_\Lambda\Lambda$, then clearly $H(M) = 0$ and 
$\epsilon M = \epsilon\Lambda = kQ$ is a projective $kQ$-module. Both properties
carry over to direct sums and direct summands, thus to all projective 
$\Lambda$-modules. 

Conversely, assume that $H(M) = 0$ and let $N = \epsilon M.$  
The multiplication with $\epsilon$ yields an isomorphism from
$M/\epsilon M$ onto $\epsilon M$, thus we can identify $M/epsilon M$ with $N$. 
Let $P = N[\epsilon]$ and consider the canonical map $p\: P \to P/\epsilon P = N$.
We can lift this map to a map $p'\:P \to M$ which again has to be surjective (since 
$\epsilon M$ lies in the radical of $M$). Now $\dim P = 2\dim N =
\dim M/\epsilon M + \dim \epsilon M = \dim M.$ This shows that $p'$ is an isomorphism. \par
 \hfill $\square$
	\medskip

{\bf 4.4.} {\it Any $M$ in $\Cal L$ has a projective submodule $U$ such that $\epsilon M
\subseteq U$ and such that $M/U$ can be identified with $H(M).$}
	\medskip
Proof.  Let $M \in \Cal L.$ Let $M' = \{x \in M\mid \epsilon x = 0\}$ and $M'' = \epsilon M$.
Then $M'' \subseteq M'$, since $\epsilon^2 = 0$ and the multiplication with $\epsilon$ 
yields an isomorphism $M/M' \to M''$. Obviously, 
$M/M'$ is annihilated by $\epsilon,$ thus a $kQ$-module. On the other hand, 
$M''$ as a submodule of $M$ is torsionless. Thus
the module $M/M' \simeq M''$ belongs to $\Cal L \cap kQ-\mod$, and therefore is a projective
$kQ$-module. But this means that the inclusion map $M'/M'' \to M/M''$ is a split monomorphism
(of $\Lambda$-modules), 
since these are $kQ$-modules and the cokernel is a projective $kQ$-module. In this way, we
obtain a $\Lambda$-submodule $U$ of $M$ such that $U\cap M' = M''$ and $U+M' = M.$
Note that $M/U \simeq M'/M'' = H(M)$.

It remains to see that $U$ is a projective 
$\Lambda$-module. In order to use 4.3, we have to show that $H(U) = 0$ and that $\epsilon U$
is a projective $kQ$-module.  Let $U' = \{x\in U\mid \epsilon x = 0\}$, then
$U' = U\cap M'$, thus $U' = M''.$  Let us show that  $\epsilon U = M''$. On the one hand 
$\epsilon U \subseteq \epsilon M = M''.$ 
On the other hand, any element in $M'' = \epsilon M$ has the form
$\epsilon x $ with $x \in M = U+M'$, thus $x = u+x'$ with $u\in U$ and $x' \in M'$, thus
$\epsilon x = \epsilon(u+x')= \epsilon u$ belongs to $\epsilon U.$ It follows
that $H(U) = U'/\epsilon U = M''/M'' = 0.$ Also, we know that $\epsilon U = M''$ is 
a projective $kQ$-module.  \par\hfill $\square$ 
	\bigskip
We see: If $M$ belongs to $\Cal L$, then there is an exact sequence
$$
   0 @>>> U @>>> M @>g>> H(M) @>>> 0
$$
with $U$ a projective $\Lambda$-module. This means that we have obtained a
right $\Cal L$-approximation $g$ of $H(M)$. But note  
 that we do cannot obtain a canonical such map,
since the submodule $U$ is not uniquely determined. Also observe that this is a minimal right
$\Cal L$-approximation if and only if $M$ has no non-zero projective direct summand.
	\medskip
{\bf 4.5.} {\it If $M$ is in $\Cal L$, indecomposable and not projective, then $H(M)$ is
indecomposable and $\eta H(M)$
is isomorphic to $M$.}
	\medskip
	
Proof.  If there is a direct decomposition 
$H(M) = N_1\oplus N_2$, let $M_i \to N_i$ be a minimal right $\Cal L$-approximation
of $N_i$, for $i = 1,2$. Then $M_1\oplus M_2 \to N_1\oplus N_2$ is a minimal right 
$\Cal L$-approximation. The uniqueness of minimal right $\Cal L$-approximations yields
that $M$ is isomorphic to $M_1\oplus M_2$, thus one of the $M_i$ is zero, say $M_2 = 0.$
But then also $N_2 = 0.$ 

Since $M \to H(M)$ is a minimal right $\Cal L$-approximation, we see that $\eta H(M)$
and $M$ are isomorphic.  \hfill $\square$
	\bigskip
In 4.4 we have seen that for any module $M\in \Cal L$ there is a projective
submodule $U$ of $M$ such that $\epsilon M \subseteq U$ and such that $M/U$
can be identified with $H(M)$. Let us show that 
for any projective submodule $U'$ with $\epsilon M \subseteq U'$,  
the projection map $p\:M \to M/U'$ induces an isomorphism $H(p).$ 
	\medskip
{\bf 4.6.} {\it Let $M$ be a module in $\Cal L$. Let
$U$ be a projective submodule of $M$ such that $\epsilon M \subseteq U.$  Then
the map $H(M)\to H(M/U) = M/U$ induced by the projection is an isomorphism.}
	\medskip
Proof. First, let us show that $\epsilon M = \epsilon U.$
In order to prove this, we may consider $M$ as an $A$-module (recall that
$A = k[\epsilon]$) and $U$ as an $A$-submodule of
$M$. Since $U$ is projective as a $\Lambda$-module, $U$ is also projective as an $A$-module
(since $kQ[\epsilon]$ is projective as an $A$-module). But $A$ is self-injective, thus any
projective $A$-module is also injective as an $A$-module. This shows that the embedding 
$U \to M$ splits as an embedding of $A$-modules. Thus there is an $A$-submodule $U'$ of $M$
such that $M = U\oplus U'.$ But $U'$ is isomorphic to $M/U$ as an $A$-module, thus annihilated
by $\epsilon.$ This shows that $\epsilon M = \epsilon U \oplus \epsilon U' = \epsilon U.$

Let $M' = \{x\in M\mid \epsilon x = 0\}.$ Then 
$M'\cap U =  \{x\in U\mid \epsilon x = 0\} = \epsilon U = \epsilon M.$

Also, $ M' + U = M.$
For the proof, consider an element $x\in M$. Since  $\epsilon M = \epsilon U,$
there is $u\in U$ such that $\epsilon x = \epsilon u$, thus $\epsilon(x-u) = 0$.
This shows that $x-u\in M'$ and therefore $x = (x-u)+u\in M'+U.$

There is the canonical (Noether-) isomorphism 
$$
 M'/(M'\cap U) \to (M'+U)/U.
$$
It yields the required isomorphism
$$
 H(M) = M'/\epsilon M = M'/(M'\cap U) \simeq (M'+U)/U = M/U = H(M/U).
$$
 \hfill $\square$
	\medskip
If we start with a module $N\in \mod kQ$ and form $\eta N,$ then there is given an exact
sequence
$$
 0 @>>> U @>u>> \eta N @>g>> N @>>> 0,
$$
such that $U$ is projective and $\eta N$ belongs to $\Cal L$.
thus we can apply 4.6 in order to see
that the map $H(\eta N)\to H(N) = N$ induced by $g$ is an isomorphism. This shows:
	\medskip
{\bf 4.7.} {\it If $N\in \mod kQ,$ then $H(\eta N)$ is 
isomorphic to $N$.} 
 \hfill $\square$
	\bigskip
It remains to study the minimal right $\Cal L$-approximation for $N\in \mod kQ$.
Start with a minimal projective $kQ$-resolution of $N$, say
$$
 0 @>>> \Omega_0 N            @>>> P_0N @>>> N @>>> 0 
$$
and embed $\Omega_0 N$ into
$(\Omega_0N)[\epsilon]$. Since $\Omega_0N$ is a projective $kQ$-module, we know that
$(\Omega_0N)[\epsilon]$ is a projective $\Lambda$-module. 
Forming the pushout of the given embeddings of $\Omega_0N$ into
$P_0N$ and $(\Omega_0N)[\epsilon]$, 
we obtain the following commutative diagram 
of $\Lambda$-modules with exact rows and columns
$$
\CD
   @. 0 @.  0  \cr 
@.      @VVV                  @VVV        \cr
0 @>>> \Omega_0 N            @>>> P_0N @>>> N @>>> 0 \cr
@.      @VVV                  @VVV      @|    \cr
0 @>>> (\Omega_0N)[\epsilon] @>>> X   @>g>> N @>>> 0 \cr
@.      @VVV                  @VVV        \cr
   @>>> \Omega_0 N @=   \Omega_0 N  \cr 
@.      @VVV                  @VVV        \cr
   @. 0 @.  0  
\endCD
$$

{\bf 4.8.} {\it The module $X$ belongs to $\Cal L$ and has no non-zero 
projective direct summand. As a consequence, the map $g\:X \to N$ is a 
minimal right $\Cal L$-approximation of $N$.}
	\medskip
Thus, we can identify $X$ with $\eta N.$
	\medskip
Proof.  The vertical sequence in the middle shows that
$X$ is an extension of $P_0N$ and $\Omega_0 N$.
Both $P_0N$ and $\Omega_0 N$ belong to $\Cal L$ and $\Cal L$ is closed under extensions, thus
$X$ belongs to $\Cal L$.
On the other hand, we have already mentioned that $(\Omega_0N)[\epsilon]$ is a projective
$\Lambda$-module, thus $g$ is a right $\Cal L$-approximation.

Let us show that any projective direct summand $P$
of $X$ is zero. According to 4.6, we know that $H(X) = N.$
Now assume that there is a direct decomposition of $\Lambda$-modules $X = M\oplus P$,
with $P$ projective. Then $H(M) = H(X)$. Let 
$M' = \{x\in M\mid \epsilon x = 0\}$ and $M'' = \epsilon M.$ According to 4.2 the following
exact sequence
$$
 0 \to M'' \to M' \to H(M) \to 0
$$
is a projective $kQ$-resolution of $H(M)$. Since we are using a minimal
projective $kQ$-resolution of $N = H(X) = H(M)$, we see that there is a projective
$kQ$-module $C$ with $C \oplus \Omega_0 N = M''$ and $C\oplus P_0N = M'.$
It remains to compare the dimensions: The equality
$$
\align
    \dim M &= 2\dim M'' + \dim N \cr
           &= 2\dim C + 2 \dim \Omega_0N + \dim N \cr
           &= 2\dim C + \dim X \cr
           &= 2\dim C + \dim M + \dim P 
\endalign
$$
shows that $P$ (and $C$) have to be zero (the first line uses that $H(M) = N$).

This implies that $g$ is minimal, since otherwise there is a non-trivial direct
decomposition $X = X'\oplus X''$ 
with say $X''$ contained in the kernel of $g$. But then $X''$ is isomorphic to a direct
summand of the kernel of $g$, thus projective.  \hfill $\square$
	\medskip
{\bf 4.9. Proof of Theorem 2.} Let $M$ be an indecomposable module in $\Cal L\setminus \Cal P.$
We attach to it $H(M)$, this is a $kQ$-module. 
According to 4.5, we know that $H(M)$ is indecomposable and also that $\eta H(M) = M.$

Conversely, let us start with an indecomposable $kQ$-module $N$. By construction,
$\eta N$ belongs to $\Cal L.$ We know by 4.8 that
$\eta N$ has no non-zero projective direct summands. And we know by 4.7 that $H(\eta N) = N$.
We still have to show that $\eta N$ is indecomposable. Thus assume that 
there is a non-trivial direct decomposition $\eta N = M_1\oplus M_2$ of $\Lambda$-modules.
Both modules $M_1, M_2$ belong to $\Cal L\setminus \Cal P$, thus according to
the first part of the proof, 
$H(M_1)$ and $H(M_2)$ are non-zero. As a consequence,  
$N = H(\eta N) = H(M_1)\oplus H(M_2)$ is 
a non-trivial direct decomposition. This contradicts the assumption that $N$ is
indecomposable.

It remains to be shown that $H\: \Cal L \to \mod kQ$ is full. Let $M_1,M_2$ be modules
in $\Cal L$. We have to show that $H$ yields a surjection $\Hom_\Lambda(M_1,M_2) \to 
\Hom_{kQ}(H(M_1),H(M_2)).$ Of course, we can assume that $M_1$ and $M_2$ both are
indecomposable. If one of the modules $M_i$ is projective, then $H(M_i) = 0$ and
nothing has to be shown. Thus we can assume that both modules 
belong to $\Cal L\setminus\Cal P.$
It follows that there are $\Cal L$-approximations $g_i\:M_i \to H(M_i)$ for $i=1,2$. 
But then any
homomorphism $f\:H(M_1) \to H(M_2)$ can be lifted to a map $\widetilde f\:M_1 \to M_2$
with $g_2\widetilde f = fg_1$ and $H(\widetilde f) = f$.  \hfill $\square$
 \hfill $\square$
	\bigskip\bigskip

{\bf Illustration.} The picture to have in mind when dealing with $\eta N$ for $N\in \mod kQ$
is the following:
$$
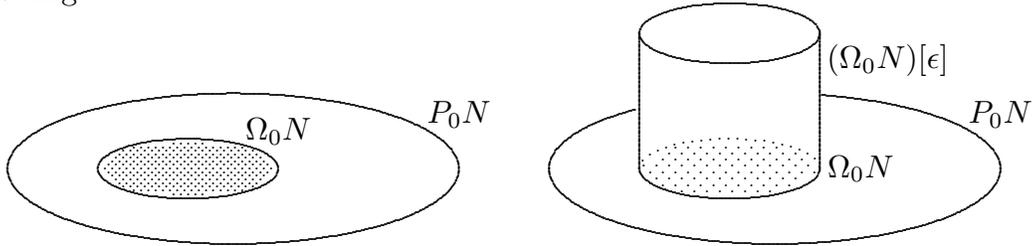
{\beginpicture
\setcoordinatesystem units <1cm,.6cm>
\put{\beginpicture
\multiput{} at 0 -.5  5 3.5 /
\ellipticalarc axes ratio 3:1 360 degrees from 0 1 center at 2.5 1  
\ellipticalarc axes ratio 3:1 360 degrees from 1 1 center at 2 1 
\put{$P_0N$} at 5 2 
\put{$\Omega_0N$} at 3 1.7 
\setshadegrid span <.4mm>
\vshade 1 1 1   <z,z,,>  1.3 .6 1.4 <z,z,,>  2 .5 1.5 <z,z,,> 2.7 .6 1.4 <z,z,,>  3 1 1 /
\endpicture} at 0 0
\put{\beginpicture
\multiput{} at 0 -.5  5 3.5 /
\ellipticalarc axes ratio 3:1 180 degrees from 0 1 center at 2.5 1  
\ellipticalarc axes ratio 3:1 -52 degrees from 0 1 center at 2.5 1  
\ellipticalarc axes ratio 3:1 78 degrees from 5 1 center at 2.5 1  
\plot 1 1  1 3.5 /
\plot 3 1  3 3.5 /
\ellipticalarc axes ratio 3:1 180 degrees from 1 1 center at 2 1 
\ellipticalarc axes ratio 3:1 360 degrees from 1 3.5 center at 2 3.5 
\setdots <1mm>
\ellipticalarc axes ratio 3:1 -180 degrees from 1 1 center at 2 1 

\put{$P_0N$} at 5 2 
\put{$\Omega_0N$} at 3.45 1 
\put{$(\Omega_0N)[\epsilon]$} at 3.75 3 
\setshadegrid span <.7mm>
\vshade 1 1 1   <z,z,,>  1.3 .6 1.4 <z,z,,>  2 .5 1.5 <z,z,,> 2.7 .6 1.4 <z,z,,>  3 1 1 /
\endpicture} at 6 0
\endpicture}
$$
On the left, we show $\Omega_0N$ as a submodule of $P_0N$, note that $N$ is obtained
from $P_0N$ by factoring out $\Omega_0N$. The right pictures depicts $\eta N$ as an amalgamation
of $P_0N$ and $(\Omega_0N)[\epsilon]$ along $\Omega_0N$. Note: if we write $M = \eta N$,
and set $M' = \{x\in M\mid \epsilon x = 0\}$ and $M'' = \epsilon M$, then $M' = P_0N$
and $M'' = \Omega_0N.$ 

The right picture looks similar to the usual depiction of a mapping cylinder in topology,
but better it should be compared with a mapping cone. After all, what here looks like a
cylinder is a projective $\Lambda$-module, thus something that one should consider as 
a contractible object. 

Of course, the same kind of pictures can be used also to depict the projective 
$\Lambda$-modules, 
stressing that they are of the form $P_0[\epsilon],$ with $P_0$ a projective
$kQ$-module.
$$
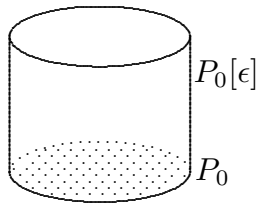
{\beginpicture
\setcoordinatesystem units <1cm,.6cm>
\put{\beginpicture
\multiput{} at 0 .5  5 3.8 /
\plot 1 1  1 3.5 /
\plot 3 1  3 3.5 /
\ellipticalarc axes ratio 3:1 180 degrees from 1 1 center at 2 1 
\ellipticalarc axes ratio 3:1 360 degrees from 1 3.5 center at 2 3.5 
\setdots <1mm>
\ellipticalarc axes ratio 3:1 -180 degrees from 1 1 center at 2 1 

\put{$P_0$} at 3.25 1 
\put{$P_0[\epsilon]$} at 3.4 2.8 
\setshadegrid span <.7mm>
\vshade 1 1 1   <z,z,,>  1.3 .6 1.4 <z,z,,>  2 .5 1.5 <z,z,,> 2.7 .6 1.4 <z,z,,>  3 1 1 /
\endpicture} at 6 0
\endpicture}
$$

{\bf Remark.} If one starts with an indecomposable non-projective $kQ$-module $N$ and 
compares the $\Lambda$-modules $N$ and $\eta N$, the decisive difference is
that $\eta N$ has $\Omega_0N$ as a factor module, and this is a non-zero projective
$kQ$-module. Thus $\Hom_\Lambda(\eta N,kQ) \neq 0$,
whereas, of course, $\Hom_\Lambda(N,kQ) = 0.$
	\bigskip
Let us insert an interesting property of the torsionless $\Lambda$-modules.
Recall that a module $M$ is said to be {\it strongly Gorenstein-projective} [BM]
provided there exists an exact sequence
$$
 0 @>>> M @>u>> P @>p>> M @>>> 0
$$
such that $u$ is a left $\Lambda$-approximation. This means, that $M$
is the kernel of a map $f\:P \to P$, where
$$
 \cdots @>>> P @>f>> P @>f>> P @>>> \cdots
$$
is a complete projective resolution (namely, $f = up$). 
	\medskip
{\bf 4.10. Proposition.} {\it Any
torsionless $kQ[\epsilon]$-module is strongly Gorenstein-projective.}
	\medskip
Proof: Let  $M$ be a torsionless $\Lambda$-module,
where $\Lambda = kQ[\epsilon]$.
It is sufficient to construct an exact sequence of the form
$$
 0 @>>> M @>u>> P @>p>> M @>>> 0
$$
with $P$ projective. Namely, since $\Lambda$ is $1$-Gorenstein, any
injective map $M \to P$ is a left $\Lambda$-approximation. Alternatively, we
construct an endomorphism $f\:P \to P$ of a projective $\Lambda$-module $P$
with $\Im(f) = \Ker(f) = M.$ 
Of course, we can
assume that $M$ is indecomposable, but also that $M$ is not projective.
Thus $M = \eta N$ for some $kQ$-module $N$. 

The construction of $\eta N$ starts with a projective presentation
$$
 0 @>>> \Omega_0N @>u>> P_0N @>>> N @>>> 0,
$$
and the map $u$ gives rise to a map $u[\epsilon]\: (\Omega_0N)[\epsilon] \to
(P_0N)[\epsilon]$, thus we may consider the map
$$
 f = \bmatrix     -\epsilon & 0 \cr
             u[\epsilon] & \epsilon
 \endbmatrix\: P \longrightarrow P,\qquad\text{where}\quad P =
(\Omega_0N)[\epsilon] \oplus (P_0N)[\epsilon].
$$
First of all, one easily verifies that $f^2 = 0$ (since $\epsilon^2 = 0$
and $\epsilon$ commutes with all maps), thus $\Im(f) \subseteq \Ker(f)$. 
We claim that the cokernel $\Cok(f)$ maps onto $\eta N$.
	\medskip
The module $\eta N$ was constructed as the following pushout
$$
\CD 
 \Omega_0N @>u>> P_0N \cr
 @V\mu VV            @VV\mu'V \cr
 (\Omega_0N)[\epsilon] @>u'>> \eta N 
\endCD
$$
Let us denote by $\pi\:(\Omega_0N)[\epsilon] \to \Omega_0$ and 
$\pi'\:(P_0N)[\epsilon] \to P_0N$
the canonical projection maps (thus $\mu\pi$ is the multiplication by
$\epsilon$ on $(\Omega_0N)[\epsilon]$ and $\pi'\cdot u[\epsilon] =
u\pi$). Let us consider
the map 
$$
 g = \bmatrix u'& \mu'\pi\endbmatrix\: 
(\Omega_0N)[\epsilon] \oplus (P_0N)[\epsilon] \longrightarrow \eta N.
$$
The composition of $f$ and $g$ is
$$
 gf = \bmatrix u'& \mu'\pi\endbmatrix \bmatrix     -\epsilon & 0 \cr
             u[\epsilon] & \epsilon
 \endbmatrix 
    = \bmatrix -u'\epsilon+\mu'\pi u[\epsilon]& \mu'\pi\epsilon\endbmatrix,
$$
but this is the zero map since $u'\epsilon = u'\mu\pi = \mu'u\pi =
\mu'\pi'u[\epsilon]$ and $\pi'\epsilon = 0$.
Therefore $g$ factors through the cokernel of $f$, thus $g$ is of the
form $P \to \Cok(f) @>\zeta >>\eta N$ for some map $\zeta\:\Cok(f) \to
\eta N.$ Since $[u',\mu']$ is surjective, also $g$ is surjective, 
and therefore $\zeta$ is surjective.
	\medskip
It remains to look at the length of $\eta N$. 
For any $\Lambda$-module $M$, let $|M|$ denote its length (and note that the 
simple $\Lambda$-modules are the simple $kQ$-modules). 
We have obtained a surjective map $\zeta\:\Cok(f) \to \eta N$, thus,
in particular, $|\eta N| \le |\Cok(f)|$. Since $\eta N$ is an extension
of $P_0N$ by $\Omega_0N$, we see that $|\eta N| = |P_0N| + |\Omega_0N|$
and therefore 
$$
 |P|  = |(P_0N)[\epsilon]| + |(\Omega_0N)[\epsilon]| =  2|\eta N|.
$$
Finally, we use that $f^2 = 0$, therefore $2|\Cok(f)| \le |P|.$
Combining the inequalities, we see that
$$
 |P| = 2|\eta N| \le 2|\Cok(f)| \le |P|.
$$
Consequently,  we must have both $|\eta N| = |\Cok(f)|$ and
$2|\Cok(f)| = |P|$. The first equality means that $\zeta$ is an
isomorphism, thus the cokernel of $f$ is isomorphic to $\eta N$.
The second equality means that $\Im(f) = \Ker(f).$  Altogether, we 
conclude that 
$f\:P \to P$ is an endomorphism of the projective $\Lambda$-module $P$
with $\Im(f) = \Ker(f) = M.$ \hfill $\square$

	\bigskip\bigskip
{\bf  5. The ghost sequences.} 	
	\medskip
Let us construct the Auslander-Reiten sequences in $\Cal L$ which end in an indecomposable
projective $kQ$-module $P_0(y)$. Note that there is the corresponding indecomposable
projective $\Lambda$-module $P(y) = (P_0(y))[\epsilon] = P_0(y)\otimes_k A$. We may identify
its submodule $P_0(y)\otimes k\epsilon$ with $P_0(y)$ and we have
$P(y)/P_0(y) \simeq P_0(y)$.

We will need the 
Auslander-Reiten translation $\tau$
in $\mod \Lambda$. We denote  by $\tau_{\Cal L}$ the  Auslander-Reiten translation 
in $\Cal L.$ 
	\medskip
{\bf 5.1. Lemma.} {\it For any vertex $y$ of $Q$, we have 
$\tau P_0(y) =  I_0(y)$ and  $\tau_{\Cal L} P_0(y) = \eta I_0(y).$}
	\medskip
Proof. In order to see that $\tau P_0(y) = I_0(y)$, one notes that
$P_0(y)$ is the cokernel of the multiplication map $\epsilon\:P(y) \to P(y)$,
thus 
$\tau P_0(y)$ has to be the kernel of the multiplication map $\epsilon\:I(y) \to I(y),$
and this is $I_0(y)$.
But for any indecomposable module $M$ in $\Cal L\setminus \Cal P$ the
Auslander-Reiten translate $\tau_{\Cal L} M$ of $M$ in $\Cal L$ 
is a non-zero direct summand of $\eta\tau M.$ Since $\eta I_0(y)$ is indecomposable,
we conclude that  $\tau_{\Cal L} P_0(y) = \eta I_0(y).$   \hfill $\square$
	\medskip
We will distinguish whether the vertex $y$ of $Q$ is a source or not.
First, we consider the case where $y$ is a source.
	\medskip
{\bf 5.2. Lemma.} {\it Let $y$ be a source of $Q$. Then the 
Auslander-Reiten sequence in $\Cal L$ 
ending in $P_0(y)$ is of the form}
$$
 0 @>>> \eta I_0(y) @>>> P(y) \oplus \rad P_0(y) @>>> P_0(y) @>>> 0.
$$
	\medskip
Proof.  We denote by $\iota\:\rad P_0(y) \to P_0(y)$ the inclusion map
and by $\epsilon'\:P(y) \to P_0(y)$ the surjective 
map induced by the multiplication with $\epsilon.$
We claim that the map $\bmatrix \epsilon' & \iota\endbmatrix\:P(y) \oplus \rad P_0(y) @>>> P_0(y)$
is minimal right almost split. 

Let $L$ be an indecomposable module in $\Cal L$
and consider a map $f\:L \to P_0(y).$ Now either $f$ maps into the radical $\rad P_0(y)$ of
$P_0(y)$, then we have a lifting. Or else it maps onto $P_0(y)$. But then there is a
map $f'\: P(y) \to L$ such that $ff' = \epsilon'$.
If the image of $f$ is annihilated by $\epsilon$, then $f'$ factors through
$\epsilon'$, say $f' = f''\epsilon'$ and then $ff''\epsilon' = ff' = \epsilon'$ 
and therefore $ff'' = 1$, thus $f$ is a split
epimorphism. Otherwise we use an embedding of $L$ into a projective $\Lambda$-module
in order to see that $f'$ has to be an isomorphism (here we use that $y$ is a source).
But then $f = \epsilon'(f')^{-1}$ is the required factorization.

It remains to show that the map is minimal. Now $P(y)$ cannot be deleted, since the
map $\epsilon'$ is surjective, whereas the map $\iota$ is not surjective. 
Also, if $N$ is an indecomposable direct summand of $\rad P_0(y)$, say $\rad P_0(y) = N\oplus 
N'$, then the inclusion map $N \to P_0(y)$ cannot be factor through the restriction of
$\bmatrix \epsilon' & \iota\endbmatrix$ to $P(y)\oplus N'$,
since the composition of any map $N \to P(y)$ with $\epsilon'$ is zero.
 \hfill $\square$
	\bigskip
{\bf 5.3. Lemma.} {\it Let $y$ be not a source of $Q$. Then the
 Auslander-Reiten sequence in $\Cal L$ 
ending in $P_0(y)$ is of the form}
$$
 0 @>>> \eta I_0(y) @>>> \eta (I_0(y)/S(y)) \oplus \rad P_0(y) @>>> P_0(y) @>>> 0.
$$
	\medskip
Proof. First, let us construct an Auslander-Reiten sequence ending in $P_0(y)$ in $\mod\Lambda$.
Since the radical of the endomorphism ring of $P_0(y)$ is zero, it is sufficient to construct
an arbitrary non-split exact sequence with end terms $I_0(y)$ and $P_0(y)$. We start with the non-split exact sequence 
$$
 0 \to S(y) \to S(y)[\epsilon]  \to S(y) \to 0,
$$
and form the induced exact sequences first with respect to the inclusion $S(y) \to I_0(y)$
and then with respect to the  projection $P_0(y) \to S(y)$ in order to obtain a sequence
$$
 0 \to I_0(y) \to M \to P_0(y) \to 0
$$
(actually, this is an Auslander-Reiten sequences with indecomposable middle term, 
as considered in [BR]). Note that
$M$ has a filtration $M_2 \subset M_1 \subset M$, with 
$$
   M_2 = \rad P_0(y), \quad M_1/M_2 = S(y)[\epsilon], \quad M/M_1 = I_0(y)/soc.
$$ 

Next, the corresponding Auslander-Reiten sequence in $\Cal L$ ending in $P_0(y)$ is
of the form
$$
 0 \to \eta I_0(y) \to \overline M \to P_0(y) \to 0,
$$
where $\overline M$ is an $\Cal L$-approximation of $M$. Our aim is to
determine $\overline M.$ The exact sequence
$0 \to M_2 \to M \to M/M_2 \to 0$ yields an exact sequence
$$
 0 \to M_2 \to \overline M \to \eta (M/M_2) \to 0, \tag{$*$}
$$ 
here we use that $M_2$ is a projective $kQ$-module, so that $\eta M_2 = M_2.$
	\medskip
Let us assume now that $y$ is not a source. 
	\medskip
{\bf Claim.} {\it We have $\eta (M/M_2) = \eta(M/M_1)$  
and the sequence $(*)$ splits.}
	\bigskip
The exact sequence
$$
 0 @>>> S(y)[\epsilon] @>>> M/M_2 @>q>> M/M_1 @>>> 0 
$$
yields an exact sequence
$$
 0 @>>> \eta(S(y)[\epsilon]) @>>> L  @>>> \eta(M/M_1) @>>> 0 
$$
where $L$ is a (not necessarily minimal) $\Cal L$-approximation of $M/M_2$.
Clearly, $\eta(S(y)[\epsilon]) = P(y)$,
but this implies that the sequence splits (since $\eta(M/M_1)\in \Cal L$ and 
$\Ext^1(\Cal L,\Cal P) = 0)$). Thus $L = P(y)\oplus \eta(M/M_1)$ is an
$\Cal L$-approximation of $M/M_2.$ Now $M/M_1 = I_0(y)/S(y)$ is the direct
sum of modules of the form $I_0(z)$ with $z\in Q_0$ (actually, the direct summands
are just all the modules $I_0(s(\alpha))$, where $\alpha$ is an arrow with $t(\alpha) = y$),
and, as we know, the minimal $\Cal L$-approximations of these indecomposable modules have
no non-zero projective direct summands. As a consequence, also a minimal $\Cal L$-approximation
of $M/M_1$ has no non-zero projective direct summand. Thus we may delete $P(y)$ and
see that $\eta(M/M_1)$ is an $\Cal L$-approximation of $M/M_2$. But actually, this has
to be a minimal $\Cal L$-approximation.
(Here is the proof. Let $g\:\eta(M/M_1) \to M/M_1$ be the minimal right $\Cal L$-approximation
and let $g'\: \eta(M/M_1) \to M/M_2$ (with $g = qg'$) be a lifting 
which is a $\Cal L$-approximation. 
Assume that we have a direct decomposition $\eta(M/M_1) = L \oplus L'$ such that 
$g'u$ is an $\Cal L$-approximation of $M/M_2$, where $u\:L \to \eta(M/M_1)$ is the inclusion map.
Then we claim that $gu$ is an $\Cal L$-approximation of $M/M_1.$ 
Assume that there is given a map $h\: L'' \to M/M_1$ with $L'' \in \Cal L.$ Then there is a lifting
$h'\: L'' \to \eta(M(M_1)$ such that $gh' = h.$ In this way, we get a map
$g'h'\:L'' \to M/M_2$. We use now that $g'u$ is an $\Cal L$-approximation of $M/M_2$ and
obtain $h''\:L'' \to L$ with $g'h' = g'uh''.$ But then $h = gh' = qg'h' = qg'uh'' = guh'$ shows
that $h$ can be factorized through $g$.)  
	
We have shown that $\eta (M/M_2) = \eta(M/M_1)$. 
In order to see that the sequence
$$
 0 \to M_2 \to \overline M \to \eta (M/M_1) \to 0
$$ 
splits, we use the fact that there is an exact sequence of the form
$$
 0 \to P \to \eta(M/M_2) \to M/M_2 \to 0.
$$
Of course, $\Ext^1_\Lambda(P,M_2) = 0$, since $P$ is projective. But also 
$\Ext^1_\Lambda(M/M_1,M_2) = 0,$ since the support of 
$M/M_1 = I_0(y)/S(y)$ are the proper predecessors of $y$, whereas the support of
$M_2 = \rad P_0(y)$ are the proper successors of $y$, thus these supports 
are separated by the vertex $y$. 
	\medskip 
Altogether, we see that $\overline M = M_2\oplus \eta(M/M_1) = \rad P_0(y)\oplus \eta I_0(y),$
and this is the middle term of the Auslander-Reiten sequence in $\Cal L$ ending in
$P_0(y)$, as we want to show. \par \hfill $\square$
	\bigskip\bigskip
{\bf  6. The ghost maps.} 	
	\medskip
A {\it ghost map} is by definition a homomorphism $f$ such that $H(f) = 0$. Let us
start with a characterization of the ghost maps between indecomposable $\Lambda$-modules
in $\Cal L.$
   \medskip
{\bf 6.1. Proposition.} {\it Let $X,Y$ be indecomposable $\Lambda$-modules in $\Cal L$ 
and let $f\: X \to Y$ be a homomorphism.
Let $X'$ be the kernel of the multiplication map 
$\epsilon\:X \to X$,
let $Y_P$ be a projective submodule of $Y$ such that $\epsilon Y \subseteq Y_P.$ 
Then $H(f) = 0$
if and only if $f$ can be written as the sum of two homomorphisms $f_0,f_1\:X \to Y$ such that
$f_1$ vanishes on $X'$, whereas the image of $f_0$ is contained in $Y_P$.}
      \medskip
Proof.  First let us show that maps of the form $f_0$ and $f_1$ as mentioned in the assertion
belong to the kernel of $H$. Now $H(f_1)$ is induced by the restriction of $f_1$ to $X'$,
thus if this restriction is zero, then $H(f_1) = 0.$ 
And, if the image of $f_0$ is contained in
$Y_P$, then $f_0$ factors through $Y_P$, but $Y_P$ is a projective $\Lambda$-module, thus
$H(Y_P) = 0$ and therefore $H(f_0) = 0.$

Conversely, let $f\:X \to Y$ be a homomorphism such that $H(f) = 0.$ 
Let $X'' = \epsilon X$, $Y'' = \epsilon Y$, and let $Y'$ be the kernel of the
multiplication map $\epsilon\:Y \to Y$. 
We can write $X$ and $Y$ as pushouts according to the
following diagrams (where all the maps are the canonical inclusion maps).
$$
\CD
 X'' @>u'>> X' \cr
 @Vu''VV      @VVuV \cr
 X''[\epsilon] @>>u'''> X
\endCD
\qquad \text{and} \qquad  
\CD
 Y'' @>v'>> Y' \cr
 @Vv''VV      @VVvV \cr
 Y''[\epsilon] @>>v'''> Y
\endCD
$$
Note that we may identify $Y''[\epsilon]$ with the submodule $Y_P$ of $Y$.

Now $H(f) = 0$ means that $f(X') \subseteq \epsilon Y.$ 
We denote by $f'\:X' \to Y''$ the restriction of $f$ to $X'$, 
it satisfies $vv'f'= fu.$ Also, let $f''\:X'' \to Y''$ be the restriction of $f$ to $X''$,
thus $f'' = f'u'.$ 
We can extend the map $f''\:X'' \to Y''$ to a map 
$f''[\epsilon]\:X''[\epsilon] \to Y''[\epsilon]$ such that $v''f'' = f''[\epsilon]u''$. 
Since $X$ is the pushout of the maps $u'$
and $u'',$ and $v''f'u' = v''f'' = f''[\epsilon]u''$,  
we obtain a map $h\:X \to Y''[\epsilon]$
such that $hu = v''f'$ and $hu''' = f''[\epsilon]$.
Let $f_0 = v'''h.$ Then $(f-f_0)u = fu - v'''hu = vv'f'-v'''v''f' = (vv'-v'''v'')f' = 0.$
This shows that $f_1 = f-f_0$ vanishes on $X'$. 
Altogether we see that $f = f_0+f_1$, that $f_1$ vanishes on $X_1$ and that
the image of $f_0$ is contained in the image of $h$, thus in $Y_P$.  \hfill $\square$
    \medskip
If we look at the two maps $f_0$ and $f_1$, 
$$
{\beginpicture
\setcoordinatesystem units <2.5cm,1cm>
\put{$X$} at 0 0
\put{$Y$} at 2 0
\put{$Y''[\epsilon]$} at 1 1
\put{$X/X'$} at 1 -1
\arr{0.2 0.2}{0.7 0.8}
\arr{1.3 0.8}{1.8 0.2}
\arr{1.3 -.8}{1.8 -.2}
\arr{0.2 -.2}{0.7 -.8}
\setsolid
\setquadratic
\plot .3 1  1 1.6  1.7 1 /
\arr{1.65 1.07}{1.7 1}
\plot .3 -1  1 -1.6  1.7 -1 /
\arr{1.65 -1.07}{1.7 -1}
\put{$f_0$} at 1 1.9
\put{$f_1$} at 1 -1.9
\endpicture}
$$
we see that they are of completely different nature:
Namely, $f_0$ factors through a projective $\Lambda$-module (namely $Y_P = Y''[\epsilon]$), 
thus through an object which
vanishes under $H$, whereas $f_1$ is a factorization inside the stable category 
$\underline{\Cal L}.$ The following special case is of interest: 
		 \medskip
{\bf 6.2. Corollary.} {\it Let $X$ be an indecomposable $\Lambda$-modules in $\Cal L$
and denote by  $X'$ the kernel of the multiplication map 
$\epsilon\:X \to X$.
If $Y$ is a projective $kQ$-module, and $f\:X \to Y$ a homomorphisms, then
$H(f) = 0$ if and only if $f(X') = 0.$}
      \medskip
Proof.  Since $Y$ is a projective $kQ$-module, the only projective $\Lambda$-module
contained in $Y$ is the zero module.  \hfill $\square$
	\bigskip
{\bf 6.3. Lemma.} {\it Let $M$ be indecomposable in $\Cal L \setminus \Cal P$ and let
$M'$ be the kernel of the multiplication map $\epsilon.$.
Let $Y$ be a projective $kQ$-module. Then, any homomorphism $f\:M \to Y$
vanishes on $M'$, thus is a ghost map.}
	\medskip
Proof. Let $M'' = \epsilon M.$ Since $Y$ is a $kQ$-module, we see that $M'' \subseteq \Ker(f)$.
Using the Noether isomorphism, we obtain an embedding
$$
   M'/(M'\cap \Ker(f)) \to (M'+\Ker(f))/\Ker(f) \subseteq Y.
$$
Now $Y$ is a projective $kQ$-module, thus also $M'/(M'\cap \Ker(f)$ is a projective
$kQ$-module. It is a factor module of $M'/M'' = H(M)$, thus a direct summand
of $H(M)$. But $M$ is indecomposable and
in $\Cal L \setminus \Cal P$, thus we know that $H(M)$ has no non-zero projective direct 
summands. This shows that $M'/(M'\cap \Ker(f) = 0$, thus 
$M' \subseteq \Ker(f).$  \hfill $\square$
	\medskip
{\bf 6.4.} Given a vertex $y$ of $Q$ which is not a source,
we have exhibited in 5.3 an Auslander-Reiten
sequence of $\Cal L$ ending in $P(y)$. Let us denote by $u(y)\:\rad P_0(y) \to P_0(y)$
the inclusion map and choose some map 
$c(y)\:\eta (I_0(y)/S(y)) \to P_0(y)$ such that 
$$
 [c(y),u(y)]\: \eta (I_0(y)/S(y)) \oplus \rad P_0(y) @>>> P_0(y)
$$
is a minimal right almost split map. If the vertex $y$ is a source, then we may denote by $c(y)$
the zero map $0 = \eta (I_0(y)/S(y)) \oplus \rad P_0(y) \to P_0(y).$
	\medskip
Then we have:
	\medskip
{\bf Theorem 3.} {\it The ideal of ghost maps in $\Cal L$ is generated by the
identity maps of the indecomposable projective $\Lambda$-modules as well as the 
maps $c(y)$ for the vertices $y$ of $Q$.}
	\medskip
Proof. Let $\Cal I$ be the ideal in $\Cal L$ generated by the 
the identity maps of the indecomposable projective $\Lambda$-modules as well as the
maps $c(y)$ with $y\in Q_0$. Of course, all the maps
in $\Cal I$ are ghost maps. 

(a) First, let us show that {\it all the maps $M \to Y$, where $M$ is indecomposable in $\Cal L\setminus \Cal P$ and not a projective $kQ$-module, whereas $Y$
is a projective $kQ$-module,
belong to $\Cal I.$} For the proof, we can assume that $Y$ is also indecomposable, say
$Y = P_0(y)$ for some vertex $y$.  
We use induction on $l(y)$, where $l(y)$ is the maximal length of
a path in $Q$ ending in $y$. If $(y) = 0$, then $y$ is a sink, and therefore
$\rad P_0(y) = 0.$ According to 5.3, the Auslander-Reiten sequence in $\Cal L$ ending in $y$ 
shows that the minimal right almost split map for $P_0(y)$ is of the form
$$
 c(y)\:\eta(I_0(y)/soc) \to P_0(y).
$$
Since we assume that $M$ is not projective, we can factor
$f$ through $c(y)$, but $c(y)$ belongs to $\Cal I,$ therefore $f$ belongs to $\Cal I.$

Next, assume that $l(y) > 0$. If $y$ is not a source, then
the right almost split map ending in $P_0(y)$ is of the form
$$
 g = [c(y),u(y)]\:\eta(I_0(y)/soc)\oplus \rad P_0(y) \to P_0(y),
$$
again according to 5.3.
We factor $f$ through $g$, thus we can write $f$ as a sum of maps,
where one factors through the map $c(y)$, thus belongs to
$\Cal I$, whereas the other maps
factor through an indecomposable direct summand of $\rad P_0(y)$. But all the indecomposable
direct summands of $\rad P_0(y)$ are of the form $P_0(x)$ with $l(x) < l(y)$, thus by
induction we know already that the maps $M \to P_0(x)$ belong to $\Cal I,$ 
and therefore $f\in\Cal I.$

It remains to look at the case where $y$ is a source, so that the 
right almost split map ending in $P_0(y)$ is of the form
$$
 g\:P(y) \oplus \rad P_0(y) \to P_0(y),
$$
now using 5.2.
Again, we factor $f$ through $g$, thus we can write $f$ as a sum of maps,
where one map factors through the projective module $P(y)$, and therefore belongs to
$\Cal I$, whereas the other maps
factor through a direct summand $P_0(x)$ 
of $\rad P_0(y)$. Again, we must have $l(x) < l(y)$, thus by
induction all the maps $M \to P_0(x)$ belong to $\Cal I.$ This shows again that $f$
belongs to $\Cal I.$

(b) Now let us consider arbitrary modules $X,Y$ in $\Cal L$ 
and let $f\:X \to Y$ be a ghost map. We want to show that $f$ belongs to
$\Cal I.$ We can assume that both modules $X$ and $Y$ are indecomposable. Also we can
assume that none of the modules $X,Y$ belongs to $\Cal P$. 

Let us exclude the case that $X$ is a projective $kQ$-module. In that case
$H(f) = 0$ means that $f(X) \subseteq \epsilon Y$,
but then we write $f$ as the composition of the following three maps
$$
  X @>>> X[\epsilon]  @>f[\epsilon]>> (\epsilon Y)[\epsilon] @>>> Y
$$ 
where the last map is some embedding (we know that such an embedding exists).
This shows that $f$ factors through a projective $\Lambda$-module.

(c) Thus it remains to consider the following setting: There is given a ghost map
$f\:X \to Y$, where
$X,Y$ are indecomposable modules in $\Cal L\setminus \Cal P$, and $X$ is not a projective
$kQ$-module. Let $X'$ be the  
kernel of the multiplication map $\epsilon\:X \to X$.
According to 6.1, we can write $f = f_0+f_1$ where $f_0$ factors through
a projective $\Lambda$-module and where $f_1$ vanishes on $X'$. 
Now $f_0$ belongs to $\Cal I$, thus it remains to be seen that $f_1$ is in $\Cal I.$
 Since $f_1$ vanishes on $X'$, we can
factor $f_1$ as
$$
  X @>>> X/X' @>>> Y.
$$
Now $X$ is indecomposable and not a projective $kQ$-module, whereas $X/X'$ is a projective
$kQ$-module, thus we have seen in (a) that the map $X \to X/X'$ belongs to $\Cal I.$
This shows that also $f_1$ belongs to $\Cal I$ and thus $f$ is in $\Cal I$.
 \hfill $\square$
	\bigskip
{\bf 6.5. Corollary.} {\it The ideal $I$ in $\underline{\Cal L}$ of all ghost maps
is a finitely generated ideal with $I^2 = 0.$}
	\medskip
Proof: According to Theorem 3, the ideal $I$ is generated by the residue classes
classes of the maps $c(y)$, with $y\in Q_0$, 
thus it is finitely generated. It remains to be seen that $I^2 = 0.$ Thus, let $y,y'$
be vertices and consider a composition of maps
$$
  \eta(I_0(y)/S(y)) @>>> P_0(y) @>g>> \eta(I_0(y')/S(y')) @>f>> P_0(y').
$$
Let $X = \eta(I_0(y')/S(y'))$ and $X'$ the kernel of the multiplication map
$\epsilon\:X \to X.$ Since $\epsilon$ vanishes on $P_0(y)$, it also vanishes on the image of $g$, thus the image of $g$ is contained in $X'$. On the other hand,
according to Lemma 6.3, $f$ vanishes on $X'$, thus
$fg = 0.$ \hfill $\square$

	\bigskip\medskip
{\bf 6.6.} Finally, let us describe the maps $c(y)$ in terms of
the arrows of the quiver.
	\medskip
Given a vertex $y$, we may decompose $\eta(I_0(y)/S(y))$ as the direct sum of the modules
$I_0(s(\alpha))$, where $\alpha$ runs through the arrows with $t(\alpha) = y$.
Thus the map $c(y)$ considered in 6.4 is given by maps $I_0(s(\alpha)) \to P_0(y)$,
for the various arrows $\alpha\:s(\alpha) \to y.$ Let us construct such maps 
explicitly.
	\medskip
Thus, for every arrow $\alpha\:i \to j$ in $Q$, we 
want to construct a map
$$
 c(\alpha)\:\eta I_0(i) \to P_0(j)
$$ 
which is irreducible in $\Cal L.$
	\medskip
In order to define $\eta I_0(i)$, we
start with a minimal projective $kQ$-presentation
$$
 0 \to \Omega_0 I_0(i) \to P_0I_0(i) \to I_0(i) \to 0,
$$
this yields a canonical map $\eta I_0(i) \to \Omega_0 I_0(i).$
Now consider any arrow $\alpha\:i \to j$
in $Q$. There is up to isomorphism a unique indecomposable $kQ$-module $N(\alpha)$ with
a non-split exact sequence
$$
 0 \to P_0(j) \to N(\alpha) \to I_0(i) \to 0.
$$
(Since there is an arrow $i\to j$ and $Q$ is directed, 
the supports of $P_0(j)$ and $I_0(j)$ do not intersect and
$P_0(j)_j = k$, $I_0(i)_i = k$. Thus we define $N(\alpha)$ by using for $N(\alpha)_\alpha$
the identity map $k \to k.$)

Since $P_0I_0(i)$ is projective, and $N(\alpha) \to I_0(i)$ is an epimorphism, 
we obtain the following commutative diagram:
$$
\CD
  0 @>>>  \Omega_0 I_0(i) @>>> P_0I_0(i) @>>>  I_0(i) @>>> 0 \cr
  @.       @VVc'(\alpha)V         @VVV            @| \cr
 0  @>>>  P_0(j) @>>>  N(\alpha)  @>>>  I_0(i)  @>>>  0
\endCD
$$
Note that the map $c'(\alpha)$ has to be surjective (namely, consider the induced
exact sequence using the projection map $\pi\: P_0(j) \to S(j)$; by the construction of $N(\alpha)$,
this sequence does not split, thus the composition of $c(\alpha)$ and $\pi$ cannot be
the zero map).

The required map $c(\alpha)\:\eta I_0(i) \to P(j)$ is the composition of the canonical map
$\eta I_0(i) \to \Omega_0 I_0(i)$ and $c'(\alpha)$:
$$
  \eta I_0(i) @>>> \Omega_0 I_0(i) @>c'(\alpha)>> P(j).
$$
	\medskip
One may also use a combined construction in order to deal with the various arrows
$\beta$ starting in $y$ at the same time, or also to deal with the various arrows
$\alpha$ ending in $y.$
Let $N(y)$ be obtained
by identifying the modules $I_0(y)$ and $P_0(y)$ at the vertex $y$, thus there are the 
following two exact sequences
$$
\CD
 0 @>>> P_0(y) @>>> N(y) @>>>I_0(y)/S(y) @>>> 0 \cr \cr
 0 @>>> \rad P_0(y) @>>> N(y) @>>>I_0(y) @>>> 0 
\endCD
$$
Of course, $I_0(y)/S(y)$ is  
the direct sum of the modules $I_0(s(\alpha))$ where
$\alpha$ is an arrow with $t(\alpha) = y,$ whereas 
$\rad P_0(y)$ is the direct sum of the modules $P_0(t(\beta))$ where
$\beta$ is an arrow with $s(\beta) = y.$

As above, we take a minimal projective resolution of $I_0(y)/S(y)$ or of $I_0(y)$
and obtain maps $c'(y)\:\Omega_0 (I_0(y)/S(y)) \to P_0(y)$
and $d'(y)\:\Omega_0 I_0(y) \to \rad P_0(y)$ as follows: 
$$
\CD
  0 @>>>  \Omega_0 (I_0(y)/S(y)) @>>> P_0(I_0(y)/S(y)) @>>>  I_0(y)/S(y) @>>> 0 \cr
  @.       @VVc'(y)V         @VVV            @| \cr
 0  @>>>  P_0(y) @>>>  N(y)  @>>>  I_0(y)/S(y)  @>>>  0
\endCD
$$
$$
\CD
  0 @>>>  \Omega_0 I_0(y) @>>> P_0I_0(y) @>>>  I_0(y) @>>> 0 \cr
  @.       @VVd'(y)V         @VVV            @| \cr
 0  @>>>  \rad P_0(y) @>>>  N(y)  @>>>  I_0(y)  @>>>  0
\endCD
$$

The composition of $c'(y)$ with the canonical map $\eta (I_0(y)/S(y)) \to 
\Omega_0 (I_0(y)/S(y))$ yields a map $c(y)\:\eta (I_0(y)/S(y)) \to P_0(y)$.
Similarly, we compose $d'(y)$ with the canonical map $\eta I_0(y) \to \Omega_0 I_0(y)$ 
and obtain $d(y)\:\eta I_0(y) \to \rad P_0(y)$.
These maps $c(y)$ and $d(y)$
can be used in the Auslander-Reiten sequences exhibited in 5.2 and 5.3.
	\bigskip\bigskip
{\bf  7. The position of the indecomposable projective modules.}
	\medskip
{\bf 7.1.} {\bf Lemma.} {\it Let $P(x)$ be the indecomposable projective $\Lambda$-module
corresponding to the vertex $x$. Then $H(\rad P(x)) = S(x)$.}
	\medskip
Proof. We write $P(x) = P_0(x)[\epsilon],$ thus
there is an exact sequence of the following form
$$
 0 \to (\rad P_0(x))[\epsilon] \to P_0(x)[\epsilon] \to S(x)[\epsilon] \to 0.
$$
This implies that we obtain for the radical of $P(x)$ an exact sequence of he form  
$$
 0 \to (\rad P_0(x))[\epsilon] \to \rad P(x) \to S(x) \to 0,
$$ 
therefore $H(\rad P(x)) = S(x)$.  \hfill $\square$
	\medskip
Thus, the Auslander-Reiten sequence with $P(x)$ as a direct summand of the middle term
starts with $\eta S(x).$ We distinguish whether $x$ is a source or not.
	\medskip
{\bf 7.2} If $x$ is a source, then 
$I_0(x) = S(x)$, thus we deal with the Auslander-Reiten sequence exhibited in Lemma 5.2:
$$
{\beginpicture
\setcoordinatesystem units <2cm,.8cm>
\put{$\eta S(x)$} at 0 0
\put{$P(x)$} at 1 1
\put{$\rad P_0(x)$} at 1 -1
\put{$P_0(x)$} at 2 0
\arr{0.3 0.2}{0.7 0.8}
\arr{1.3 0.8}{1.7 0.2}
\arr{1.3 -.8}{1.7 -.2}
\arr{0.3 -.2}{0.7 -.8}
\endpicture}
$$
What is of interest and should be remembered
is the fact that in this case the irreducible map starting in $P(x)$ is surjective with
target  $P_0(x)$.
	\medskip
It remains to consider the case that $x$ is not a source. Let us denote by $\tau_0$
the Auslander-Reiten translation for $\mod kQ$. Since $x$ is not a source, $S(x)$
is not an injective $kQ$-module, thus $\tau^{-1}_0S(x)$ is defined. 
	 \medskip
{\bf 7.3.} {\it If $x$ is not a source, then the Auslander-Reiten sequence 
with $P(x)$ as a direct summand of the middle term starts with $\eta S(x)$ and
ends in $\eta \tau^{-1}_0S(x).$}  \hfill $\square$
	\bigskip\bigskip
{\bf  8. Examples}
	\medskip
{\bf 8.1. First example.}
As a first example, we take a bipartite quiver $Q$ (this means that all the vertices are
sinks or sources). If $\alpha\:i \to j$ is an arrow in $Q$, then $I_0(i) = S(i)$ 
and $P(j) = S(j)$ are simple $kQ$-modules. 

We
consider the following quiver $Q$ exhibited below on the left. 
On the right side, we present the decisive parts of the
preinjective component and the preprojective component of $\Gamma(kQ)$, 
already as subquivers of the translation quiver $\Bbb Z(Q^{\op})$ 
separated only by some arrows
(they are marked by dotted diagonal lines):
$$
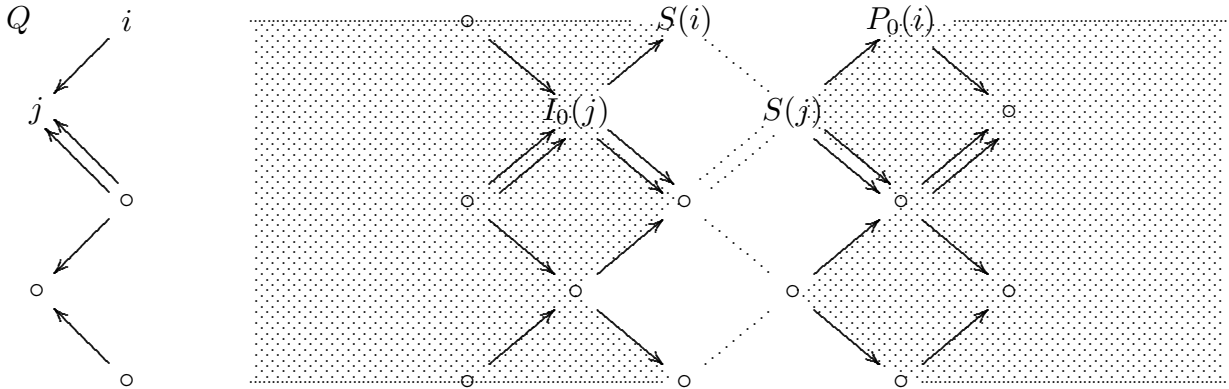
{\beginpicture
\setcoordinatesystem units <1cm,1cm>
\put{\beginpicture
\multiput{} at 0 0  1 4 /
\put{$i$} at 1 4
\put{$j$} at 0 3
\multiput{$\circ$} at 1 0  0 1  1 2 /
\arr{0.8 0.2}{0.2 0.8}
\arr{0.8 1.8}{0.2 1.2}
\arr{0.9 2.2}{0.2 2.9}
\arr{0.8 2.1}{0.1 2.8}
\arr{0.8 3.8}{0.2 3.2}
\put{$Q$} at -.2 4
\endpicture} at 0 0
\put{\beginpicture
\setcoordinatesystem units <1.2cm,1cm>
\multiput{} at 0 0  8 4 /
\put{$S(i)$} at 3 4 
\put{$I_0(j)$} at 2 3
\multiput{$\circ$} at 3 2  2 1  3 0   5 2  4 1  5 0 / 
\multiput{$\circ$} at 1 0  1 2  1 4 / 
\multiput{$\circ$} at 6 1  6 3 / 
\arr{1.2 0.2}{1.8 0.8}
\arr{1.2 1.8}{1.8 1.2}
\arr{1.3 2.1}{1.9 2.7}
\arr{1.2 2.2}{1.8 2.8}
\arr{1.2 3.8}{1.8 3.2}

\arr{2.2 0.8}{2.8 0.2}
\arr{2.2 1.2}{2.8 1.8}
\arr{2.3 2.8}{2.9 2.2}
\arr{2.2 2.7}{2.8 2.1}
\arr{2.3 3.3}{2.8 3.8}

\put{$S(j)$} at 4 3
\put{$P_0(i)$} at 5 4
\arr{4.2 0.8}{4.8 0.2}
\arr{4.2 1.2}{4.8 1.8}
\arr{4.3 2.8}{4.9 2.2}
\arr{4.2 2.7}{4.8 2.1}
\arr{4.3 3.3}{4.8 3.8}

\arr{5.2 0.2}{5.8 0.8}
\arr{5.2 1.8}{5.8 1.2}
\arr{5.3 2.1}{5.9 2.7}
\arr{5.2 2.2}{5.8 2.8}
\arr{5.3 3.7}{5.8 3.2}

\setdots <.5mm>
\plot -1 0  2.8 0 /
\plot -1 4  2.5 4 /

\plot 5.2 0  8 0 /
\plot 5.5 4  8 4 /

\setdots <1mm>
\plot 3.2 0.2  3.8 0.8 /

\plot 3.15 2.25  3.75 2.85 /
\plot 3.25 2.15  3.85 2.75 /

\plot 3.2 1.8  3.8 1.2 /
\plot 3.2 3.8  3.8 3.2 /

\setshadegrid span <.6mm>
\hshade 0 -1 3   <,,,z> 1 -1 2 <,,z,z> 2 -1 3 <,,z,z> 3 -1 2 <,,z,> 4 -1 3 /
\hshade 0 5 8   <,,,z> 1 4 8 <,,z,z> 2 5 8 <,,z,z> 3 4 8 <,,z,> 4 5 8 /

\endpicture} at 8  0
\endpicture}
$$

If we consider now $\underline{\Cal L}$, we have to replace any $kQ$-module $N$ by 
$\eta N$, and we have to add an arrow $\eta S(i) \to S(j)$ for any arrow $i \to j$
in the quiver $Q$. These new arrows represent a $k$-basis of $\Hom(\eta S(i), S(j))$.
Note that the new
arrows represent ghost maps. Here is this part of $\Gamma(\underline{\Cal L}).$
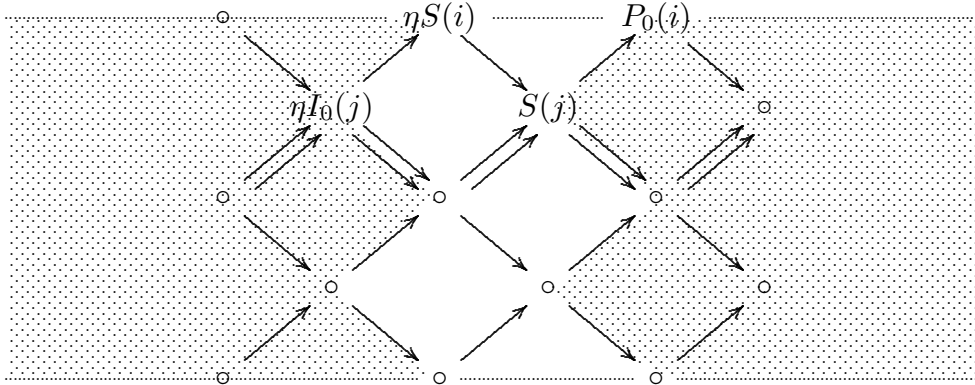
$$
{\beginpicture
\setcoordinatesystem units <1cm,1cm>
\put{\beginpicture
\setcoordinatesystem units <1.2cm,1cm>
\multiput{} at 0 0  8 4 /
\put{$\eta S(i)$} at 3 4 
\put{$\eta I_0(j)$} at 2 3
\multiput{$\circ$} at 3 2  2 1  3 0   5 2  4 1  5 0 / 
\multiput{$\circ$} at 1 0  1 2  1 4 / 
\multiput{$\circ$} at 6 1  6 3 / 
\arr{1.2 0.2}{1.8 0.8}
\arr{1.2 1.8}{1.8 1.2}
\arr{1.3 2.1}{1.9 2.7}
\arr{1.2 2.2}{1.8 2.8}
\arr{1.2 3.8}{1.8 3.2}

\arr{2.2 0.8}{2.8 0.2}
\arr{2.2 1.2}{2.8 1.8}
\arr{2.3 2.8}{2.9 2.2}
\arr{2.2 2.7}{2.8 2.1}
\arr{2.3 3.3}{2.8 3.8}

\put{$S(j)$} at 4 3
\put{$P_0(i)$} at 5 4
\arr{4.2 0.8}{4.8 0.2}
\arr{4.2 1.2}{4.8 1.8}
\arr{4.3 2.8}{4.9 2.2}
\arr{4.2 2.7}{4.8 2.1}
\arr{4.3 3.3}{4.8 3.8}

\arr{5.2 0.2}{5.8 0.8}
\arr{5.2 1.8}{5.8 1.2}
\arr{5.3 2.1}{5.9 2.7}
\arr{5.2 2.2}{5.8 2.8}
\arr{5.3 3.7}{5.8 3.2}

\arr{3.2 0.2}{3.8 0.8}
\arr{3.2 1.8}{3.8 1.2}
\arr{3.3 2.1}{3.9 2.7}
\arr{3.2 2.2}{3.8 2.8}
\arr{3.2 3.8}{3.8 3.2}

\setdots <.5mm>
\plot -1 0  2.8 0 /
\plot -1 4  2.5 4 /

\plot 3.2 0  4.8 0 /
\plot 3.5 4  4.5 4 /

\plot 5.2 0  8 0 /
\plot 5.5 4  8 4 /

\setshadegrid span <.6mm>
\hshade 0 -1 3   <,,,z> 1 -1 2 <,,z,z> 2 -1 3 <,,z,z> 3 -1 2 <,,z,> 4 -1 3 /
\hshade 0 5 8   <,,,z> 1 4 8 <,,z,z> 2 5 8 <,,z,z> 3 4 8 <,,z,> 4 5 8 /

\endpicture} at 8  0
\endpicture}
$$

Finally, let us show the corresponding part of $\Gamma(\Cal L).$ Here, the
indecomposable projective $\Lambda$-modules are added. 
$$
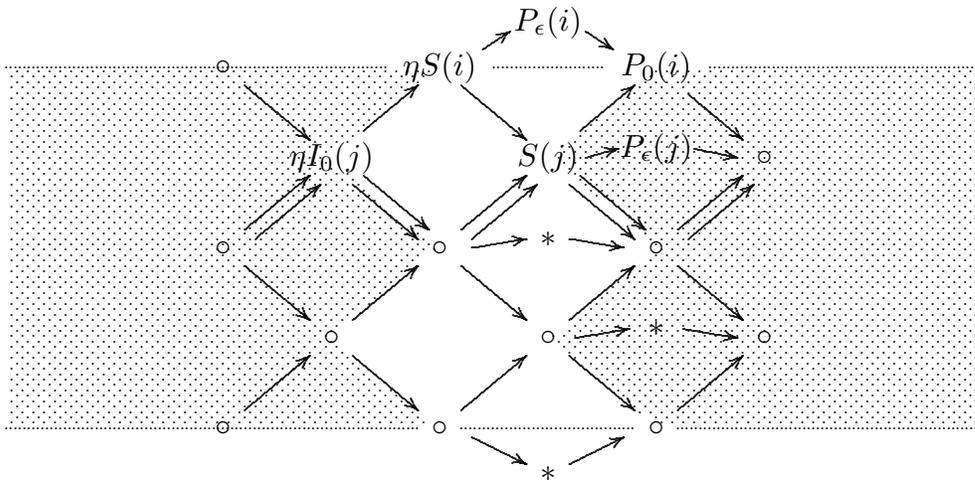
{\beginpicture
\setcoordinatesystem units <1cm,1cm>
\put{\beginpicture
\setcoordinatesystem units <1.2cm,1cm>
\multiput{} at 0 0  8 4 /
\put{$\eta S(i)$} at 3 4 
\put{$\eta I_0(j)$} at 2 3
\multiput{$\circ$} at 3 2  2 1  3 0   5 2  4 1  5 0 / 
\multiput{$\circ$} at 1 0  1 2  1 4 / 
\multiput{$\circ$} at 6 1  6 3 / 

\put{$P_\epsilon(i)$} at 4 4.5
\put{$P_\epsilon(j)$} at 5 3.1
\multiput{$*$} at 4 2.1  4 -.5  5 1.1 /
\arr{3.4 4.2}{3.65 4.4}
\arr{4.35 4.4}{4.6 4.2}

\arr{3.3 2.0}{3.8 2.1}
\arr{4.2 2.1}{4.7 2.0}

\arr{3.3 -.1}{3.8 -.4}
\arr{4.2 -.4}{4.7 -.1}

\arr{4.35 3.0}{4.65 3.1}
\arr{5.35 3.1}{5.75 3.0}

\arr{4.25 1.0}{4.75 1.1}
\arr{5.25 1.1}{5.75 1.0}

\arr{1.2 0.2}{1.8 0.8}
\arr{1.2 1.8}{1.8 1.2}
\arr{1.3 2.1}{1.9 2.7}
\arr{1.2 2.2}{1.8 2.8}
\arr{1.2 3.8}{1.8 3.2}

\arr{2.2 0.8}{2.8 0.2}
\arr{2.2 1.2}{2.8 1.8}
\arr{2.3 2.8}{2.9 2.2}
\arr{2.2 2.7}{2.8 2.1}
\arr{2.3 3.3}{2.8 3.8}

\put{$S(j)$} at 4 3
\put{$P_0(i)$} at 5 4
\arr{4.2 0.8}{4.8 0.2}
\arr{4.2 1.2}{4.8 1.8}
\arr{4.3 2.8}{4.9 2.2}
\arr{4.2 2.7}{4.8 2.1}
\arr{4.3 3.3}{4.8 3.8}

\arr{5.2 0.2}{5.8 0.8}
\arr{5.2 1.8}{5.8 1.2}
\arr{5.3 2.1}{5.9 2.7}
\arr{5.2 2.2}{5.8 2.8}
\arr{5.3 3.7}{5.8 3.2}

\arr{3.2 0.2}{3.8 0.8}
\arr{3.2 1.8}{3.8 1.2}
\arr{3.3 2.1}{3.9 2.7}
\arr{3.2 2.2}{3.8 2.8}
\arr{3.2 3.8}{3.8 3.2}

\setdots <.5mm>
\plot -1 0  2.8 0 /
\plot -1 4  2.5 4 /

\plot 3.2 0  4.8 0 /
\plot 3.5 4  4.5 4 /

\plot 5.2 0  8 0 /
\plot 5.5 4  8 4 /

\setshadegrid span <.6mm>
\hshade 0 -1 3   <,,,z> 1 -1 2 <,,z,z> 2 -1 3 <,,z,z> 3 -1 2 <,,z,> 4 -1 3 /
\hshade 0 5 8   <,,,z> 1 4 8 <,,z,z> 2 5 8 <,,z,z> 3 4 8 <,,z,> 4 5 8 /

\endpicture} at 8  0
\endpicture}
$$
	\bigskip
{\bf 8.2. Second example.} As a second example we take the quiver $Q$ of type $\Bbb A_3$
with linear orientation. First, let us show two copies of $\Gamma(\mod kQ)$ 
appropriately embedded into the translation quiver $\Bbb Z\Bbb A_3.$

$$
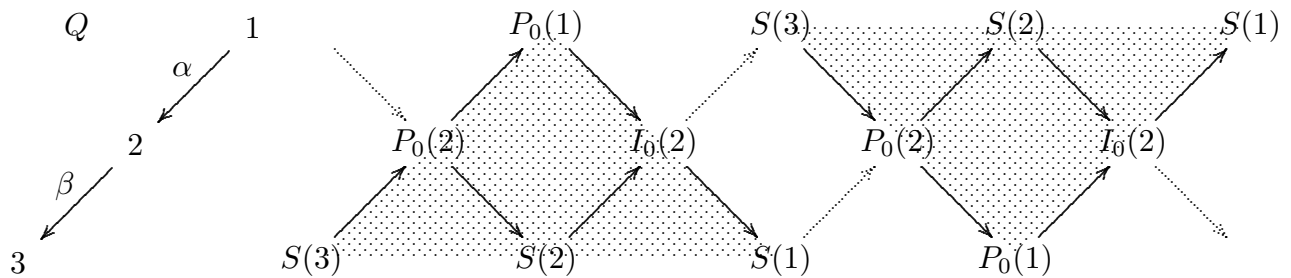
{\beginpicture
\setcoordinatesystem units <1.3cm,1.3cm>
\put{\beginpicture
\put{$3$} at 0 0
\put{$2$} at 1 1 
\put{$1$} at 2 2
\put{$\beta$} at 0.4 0.65
\put{$\alpha$} at 1.4 1.65
\arr{0.8 0.8}{0.2 0.2}
\arr{1.8 1.8}{1.2 1.2}
\put{$Q$} at 0.5 2
\endpicture} at -3.5 0

\put{\beginpicture
\put{$S(3)$} at 0 0
\put{$S(2)$} at 2 0  
\put{$S(1)$} at 4 0
\put{$P_0(2)$} at 1 1
\put{$P_0(1)$} at 2 2 
\put{$I_0(2)$} at 3 1

\put{} at 5 0
\arr{0.2 0.2}{0.8 0.8}
\arr{2.2 0.2}{2.8 0.8}
\arr{1.2 1.2}{1.8 1.8}
\arr{1.2 0.8}{1.8 0.2}
\arr{3.2 0.8}{3.8 0.2}
\arr{2.2 1.8}{2.8 1.2}

\put{$S(3)$} at 4 2
\put{$S(2)$} at 6 2  
\put{$S(1)$} at 8 2
\put{$P_0(2)$} at 5 1
\put{$P_0(1)$} at 6 0 
\put{$I_0(2)$} at 7 1

\arr{5.2 1.2}{5.8 1.8}
\arr{6.2 0.2}{6.8 0.8}
\arr{7.2 1.2}{7.8 1.8}

\arr{4.2 1.8}{4.8 1.2}
\arr{5.2 0.8}{5.8 0.2}
\arr{6.2 1.8}{6.8 1.2}

\setdots <.5mm>
\arr{0.2 1.8}{0.8 1.2}
\arr{3.2 1.2}{3.8 1.8}
\arr{4.2 0.2}{4.8 0.8}
\arr{7.2 0.8}{7.8 0.2}

\setshadegrid span <.6mm>
\hshade 0 0 4   2 2 2 /
\hshade 0 6 6   2 4 8 /
\endpicture} at  2 0

\endpicture}
$$

Now we present the corresponding torsionless $\Lambda$-modules and insert the
ghost arrows $c(\alpha)\:\eta I_0(1) \to P(2)$ and $c(\beta)\:\eta I_0(2) \to S(3).$ 
$$
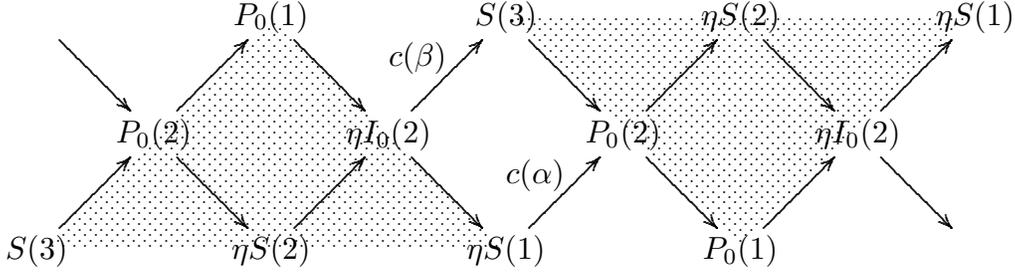
{\beginpicture
\setcoordinatesystem units <1.3cm,1.3cm>
\put{\beginpicture
\put{$S(3)$} at 0 0
\put{$\eta S(2)$} at 2 0  
\put{$\eta S(1)$} at 4 0
\put{$P_0(2)$} at 1 1
\put{$P_0(1)$} at 2 2 
\put{$\eta I_0(2)$} at 3 1

\put{$c(\beta)$} at 3.25 1.65
\put{$c(\alpha)$} at 4.25 .65

\put{} at 5 0
\arr{0.2 0.2}{0.8 0.8}
\arr{2.2 0.2}{2.8 0.8}
\arr{1.2 1.2}{1.8 1.8}
\arr{1.2 0.8}{1.8 0.2}
\arr{3.2 0.8}{3.8 0.2}
\arr{2.2 1.8}{2.8 1.2}

\put{$S(3)$} at 4 2
\put{$\eta S(2)$} at 6 2  
\put{$\eta S(1)$} at 8 2
\put{$P_0(2)$} at 5 1
\put{$P_0(1)$} at 6 0 
\put{$\eta I_0(2)$} at 7 1

\arr{5.2 1.2}{5.8 1.8}
\arr{6.2 0.2}{6.8 0.8}
\arr{7.2 1.2}{7.8 1.8}

\arr{4.2 1.8}{4.8 1.2}
\arr{5.2 0.8}{5.8 0.2}
\arr{6.2 1.8}{6.8 1.2}

\arr{0.2 1.8}{0.8 1.2}
\arr{3.2 1.2}{3.8 1.8}
\arr{4.2 0.2}{4.8 0.8}
\arr{7.2 0.8}{7.8 0.2}

\setshadegrid span <.6mm>
\hshade 0 0 4   2 2 2 /
\hshade 0 6 6   2 4 8 /
\endpicture} at  2 0

\endpicture}
$$

\noindent
The Auslander-Reiten quiver of $\Cal L$ (or better its universal covering) looks as follows: 
$$
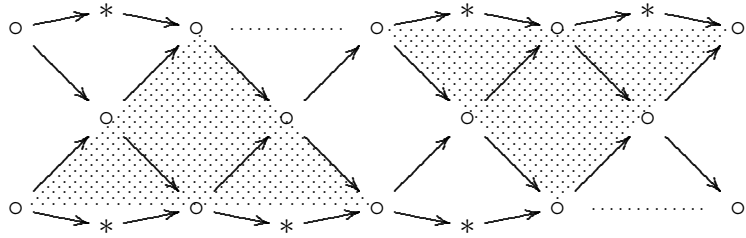
{\beginpicture
\setcoordinatesystem units <1cm,1cm>
\put{\beginpicture
\multiput{$\circ$} at 0 0  1 1  2 2  2 0  3 1  4 0 /

\multiput{$*$} at 1 -.2  3 -.2  5 -.2  
                     1 2.2   5 2.2  7 2.2  /
\arr{0.2 -.04}{0.8 -.16}
\arr{1.2 -.16}{1.8 -.04}
\arr{2.2 -.04}{2.8 -.16}
\arr{3.2 -.16}{3.8 -.04}
\arr{4.2 -.04}{4.8 -.16}
\arr{5.2 -.16}{5.8 -.04}

\arr{0.2 2.04}{0.8 2.16}
\arr{1.2 2.16}{1.8 2.04}
\arr{4.2 2.04}{4.8 2.16}
\arr{5.2 2.16}{5.8 2.04}
\arr{6.2 2.04}{6.8 2.16}
\arr{7.2 2.16}{7.8 2.04}

\setdots <1mm>
\plot 2.4 2  3.6 2 /
\plot 6.4 0  7.6 0 /
\setsolid
\put{} at 5 0
\arr{0.2 0.2}{0.8 0.8}
\arr{2.2 0.2}{2.8 0.8}
\arr{1.2 1.2}{1.8 1.8}
\arr{1.2 0.8}{1.8 0.2}
\arr{3.2 0.8}{3.8 0.2}
\arr{2.2 1.8}{2.8 1.2}

\multiput{$\circ$} at 0 2  4 2  5 1  6 0  6 2  7 1  8 2  8 0  /
\arr{5.2 1.2}{5.8 1.8}
\arr{6.2 0.2}{6.8 0.8}
\arr{7.2 1.2}{7.8 1.8}

\arr{4.2 1.8}{4.8 1.2}
\arr{5.2 0.8}{5.8 0.2}
\arr{6.2 1.8}{6.8 1.2}

\arr{0.2 1.8}{0.8 1.2}
\arr{3.2 1.2}{3.8 1.8}
\arr{4.2 0.2}{4.8 0.8}
\arr{7.2 0.8}{7.8 0.2}

\setshadegrid span <.6mm>
\hshade 0 0 4   2 2 2 /
\hshade 0 6 6   2 4 8 /
\endpicture} at  2 0

\endpicture}
$$

It may be helpful for the reader to use this example in order to write down 
all the $\Lambda$-modules as representations of $Q$
over $A$ (for example, 
the module $M = \eta Q(1) = \eta S(1)$ is written as $AAk$, since $M_3 = M_2 = A$
and $M_1 = k$). As in some illustrations before, we use dotted arrows in order to
mark the ghost maps. 
$$
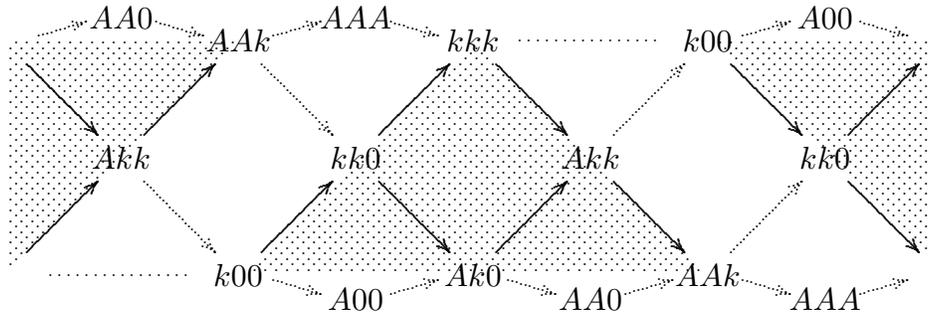
{\beginpicture
\setcoordinatesystem units <1.3cm,1.3cm>
\put{\beginpicture
\put{$AAk$} at 0 2
\put{$k00$} at 0 0
\put{$Ak0$} at 2 0  
\put{$AAk$} at 4 0
\put{$kk0$} at 1 1
\put{$kkk$} at 2 2 
\put{$Akk$} at 3 1

\arr{-1.8 1.8}{-1.2 1.2}

\arr{0.2 0.2}{0.8 0.8}
\arr{2.2 0.2}{2.8 0.8}
\arr{1.2 1.2}{1.8 1.8}
\arr{1.2 0.8}{1.8 0.2}
\arr{3.2 0.8}{3.8 0.2}
\arr{2.2 1.8}{2.8 1.2}

\put{$k00$} at 4 2
\put{$kk0$} at 5 1
\put{$Akk$} at -1 1

\arr{5.2 1.2}{5.8 1.8}

\arr{4.2 1.8}{4.8 1.2}
\arr{5.2 0.8}{5.8 0.2}
\arr{-.8 1.2}{-.2 1.8}
\arr{-1.8 .2}{-1.2 .8}

\setdots <.5mm>
\arr{-.8 0.8}{-.2 0.2}
\arr{0.2 1.8}{0.8 1.2}
\arr{3.2 1.2}{3.8 1.8}
\arr{4.2 0.2}{4.8 0.8}

\setshadegrid span <.6mm>
\hshade 0 -2 -2  2 -2 0 /
\hshade 0 0 4   2 2 2 /
\hshade 0 6 6   2 4 6 /

\put{$A00$} at 1 -.2
\put{$AA0$} at  3 -.2 
\put{$AAA$} at  5 -.2  
\put{$AAA$} at  1 2.2 
\put{$A00$} at  5 2.2 
\put{$AA0$} at  -1 2.2

\arr{0.3 -.04}{0.7 -.16}
\arr{1.3 -.16}{1.7 -.04}
\arr{2.3 -.04}{2.7 -.16}
\arr{3.3 -.16}{3.7 -.04}
\arr{4.3 -.04}{4.7 -.16}
\arr{5.3 -.16}{5.7 -.04}

\arr{0.3 2.04}{0.7 2.16}
\arr{1.3 2.16}{1.7 2.04}
\arr{4.3 2.04}{4.7 2.16}
\arr{5.3 2.16}{5.7 2.04}

\arr{-1.7 2.04}{-1.3 2.16}
\arr{-.7 2.16}{-.3 2.04}

\setdots <1mm>
\plot 2.4 2  3.6 2 /
\plot -1.6 0  -.4 0 /

\endpicture} at  2 0

\endpicture}
$$
	\bigskip\bigskip
{\bf 9. Remarks on the behavior of the homology functor $H$ on $\mod\Lambda$.}
	\medskip
We have seen that the homology functor has nice properties when restricted to the 
subcategory $\Cal L\setminus\Cal P$. For example, it maps indecomposable modules to
indecomposables, and non-isomorphic indecomposable modules to non-isomorphic ones.
Also, for a fixed homology dimension vector $\bold r$, there is either one indecomposable object
in $\Cal L\setminus \Cal P$ 
with $\bdim H(M) = \bold r$ or at least a $1$-parameter family. 
But all these assertions are only true for the restriction of $H$ to  $\Cal L\setminus\Cal P$.
In general, one cannot expect such a pleasant behavior.
	\medskip
{\bf 9.1.~Example.} Consider the quiver $Q$ of type $\Bbb A_3$ with two sources and a sink.
Then the $AQ$ module $M = k \to A \leftarrow k$ is indecomposable, but
$K(M) = k \to 0 \leftarrow k$ is decomposable.
	\medskip
{\bf 9.2.~Example.} Let $Q$ be the Kronecker quiver with two arrows from $2$ to $1$. Then
there is a $\Bbb P_1$-family of indecomposable $AQ$ modules $M$ with
$M_2 = k$ and $M_1 = A$. For all of them $H(M)$ is the simple module $S(2)$.
Thus here we have many non-isomorphic $AG$-modules with isomorphic homology modules.
	\medskip
{\bf 9.3.~Example.} Let $Q$ be the quiver of type $\Bbb D_4$ with subspace orientation, let
$1,2,3$ be the sources of $Q$  and $0$ the sink. 
Then there are precisely 2 isomorphism classes of indecomposable 
$AQ$-modules $M$ with $H(M) = S(1)\oplus S(2) \oplus S(3)$, for one of them $M_0 = A$,
for the other one, $M_0 = A^2.$ (For the proof, one has to observe that in this case
$\Lambda$ is a tubular algebra, and one has to study the corresponding root system
in detail.)  
	\bigskip\bigskip
{\bf 10. Generalization.}
	\medskip
Instead of looking at the path algebra $kQ$ of a quiver, one may start with an arbitrary
finite-dimensional $k$-algebra $H$ which is hereditary and take instead of $\Lambda =
kQ[\epsilon]$ the corresponding algebra $H[\epsilon]$. Observe that the proofs of
the theorems given here work in general. What is special in the quiver case is just the possibility
to describe the maps $c(y)$ in terms of
the arrows of the quiver, see 6.5. 
	\bigskip
  	\bigskip\bigskip

{\bf 11. References.}
	\medskip
\item{[AB]} M. Auslander, R. Buchweitz. The homology theory of Cohen-Macaulay approximations.
     Soc. Math. France. Mem 38 (1989), 5--37. 
\item{[AR1]} M. Auslander, I. Reiten. Applications of contravariantly finite subcategories.
   Adv. Math. 86 (1991), 111--152.
\item{[AR2]} M. Auslander, I. Reiten. Cohen-Macaulay and Gorenstein Artin algebras.
  In Representation theory of finite groups
   and finite-dimensional algebras, Birkh\"auser (1991), 221--245. 
\item{[BMRRT]} A. B. Buan, R. Marsh, M. Reineke,
   I. Reiten, G. Todorov. Tilting theory and cluster combinatorics. 
   Adv. Math. 204 (2006), 572--618.
\item{[B]} R. Buchweitz. Maximal Cohen-Macaulay modules and Tate-cohomology.
   Ms. (unpublished)
\item{[BM]} D. Bennis, N. Mahdou. Strongly Gorenstein projective, injective, and 
 flat modules. J. Pure Appl. Algebra 210 (2007), 437--445.

\item{[BR]} M.C.R. Butler, C.M.Ringel. Auslander-Reiten sequences with few middle terms, with applications to string algebras. Comm. Algebra 15 (1987), 145--179.
\item{[C]} X-W. Chen. Algebras with radical square zero are either self-injective
of CM-free. Proc. Amer. Math. Soc. 140 (2012), 93--98. 

\item{[EJ]} E. E. Enochs, O. M. G. Jenda. Relative homological algebra. De Gruyter Expo\. Math\. 30
   (2000).
\item{[H1]} D. Happel. Triangulated categories in the representation theory of finite
  dimensional algebras, LMS Lecture Notes Vo1.119, Cambridge University Press (1988).   
\item{[H2]} D. Happel. On Gorenstein algebras.
  In Representation theory of finite groups
   and finite-dimensional algebras, Birkh\"auser (1991) 389--404. 
\item{[Ka]} V. Kac. Infinite root systems, representations of graphs and invariant theory, 
   Invent. Math., 56 (1980), 57--92. 
\item{[Ke]} B. Keller. On triangulated orbit categories. Doc. Math. 10 (2005), 
    551--581.
\item{[KR]} B. Keller, I. Reiten. Cluster-tilted algebras are Gorenstein and stably Calabi-Yau.
  Adv. Math. 211 (2007), 123--151.
\item{[LZ]} X. H. Luo, P. Zhang. Monic representations and Gorenstein-projective modules.
  Pacific Journal of Mathematics. 264 (2013), 163--194.
\item{[O]} D. Orlov.
    Triangulated categories of singularities and D-branes in Landau-Ginzburg models.
    Proc. Steklov Inst. Math. 246 (3) (2004), 227--248.
\item{[RX]} C. M. Ringel, B.-L. Xiong. On rings with radical square zero. 
 Algebra and Discrete Mathematics. 14 (2012), 297 - 306. 

	\bigskip
{\rmk
C. M. Ringel\par
Department of Mathematics, Shanghai Jiao Tong University \par
Shanghai 200240, P. R. China, and \medskip 
King Abdulaziz University, \par
P O Box 80200, Jeddah, Saudi Arabia.\medskip  
e-mail: {\ttk ringel\@math.uni-bielefeld.de} \par
\medskip

P. Zhang\par
Department of Mathematics, Shanghai Jiao Tong University \par
Shanghai 200240, P. R. China, and \par 
e-mail: {\ttk  pzhang\@sjtu.edu.cn} \par}

\bye